\documentclass{article}
\title{Nearly K\"ahler geometry and Riemannian foliations
}
\author{Paul-Andi Nagy}
\date{\today}
\usepackage{amsfonts,amssymb}

\newtheorem{teo}{Theorem}[section]
\newtheorem{lema}{Lemma}[section]
\newtheorem{pro}{Proposition}[section]
\newtheorem{defi}{Definition}[section]
\newtheorem{rema}{Remark}[section]
\newtheorem{coro}{Corollary}[section]
\newtheorem{nr}{}[section]
\begin{document}
\maketitle
\abstract{\normalsize We consider strict and complete nearly K\"ahler manifolds with the canonical 
Hermitian connection. The holonomy representation of the canonical Hermitian connection is studied. We show that a 
strict and complete nearly K\"ahler is locally a Riemannian product of homogenous nearly K\"ahler spaces, twistor 
spaces over quaternionic K\"ahler manifolds and $6$-dimensional nearly K\"ahler manifolds. As an application we obtain structure results for totally geodesic 
Riemannian foliations admitting a compatible K\"ahler structure. Finally, we obtain a  classification result for 
the homogenous case, reducing a conjecture of Wolf and Gray to its 
$6$-dimensional form.
} 
\large
\section{Introduction}
Given an oriented Riemannian manifold $(M^n,g)$ a fundamental piece of data is encoded in its 
holonomy group, to be denoted by $Hol(M,g)$. Indeed, the condition that $Hol(M,g)$ be different from $SO(n)$ has strong geometric 
implications at the level of the geometric structure supported by $M$ (see \cite{pitman} for an account). The search for new geometric 
structures supporting (for example) Einstein metrics or satisfying usefull curvature identities motivated the 
introduction by A. Gray of the concept of weak holonomy (see \cite{Gray2} and \cite{SS}). However, in the case 
when the weak holonomy group of a manifold acts transitvely on its unit bundle it was shown \cite{Gray2} that only three groups, namely
$$ U(n), \ G_2 \ \mbox{and} \ Spin(9) $$ 
can potentially produce new geometric structures (other than those coming from Riemannian holonomy). This is 
certainly the case for weak holonomy $G_2$ (see 
\cite{Fri3}), where many homogenous examples were constructed. Also, it is now known that every $7$ dimensional $3$-sasakian manifold 
carry a metric with weak holonomy $G_2$ and, furthermore, there is no scarcity of such manifolds (see \cite{Boyer}).
Very recentely, the case of $Spin(9)$ was 
undertaken by Friedrich \cite{Fri1}. \par
In this paper we will be concerned with the study of manifolds 
with weak holonomy $U(n)$. These 
manifolds are called nearly K\"ahler and also appear as one of the sixteen classes of almost Hermitian 
manifolds \cite{Gray3}. Many properties of this class of manifolds are now known 
(see \cite{Gray5, Gray4, Gray1}). Note the $6$-dimensional 
case stands out because of the existence of Killing spinors (fact that is in this dimension equivalent to 
being nearly K\"ahler \cite{Grun1}) and because of carrying Einstein metrics. \par
Nearly K\"ahler manifolds can also be caractherized as almost Hermitian manifolds admitting a metric 
Hermitian connection whose torsion is parallel and totally skew. From this point of view they are interesting 
in theoretical physics \cite{Friedrich2}. Manifolds carrying a connection with totally skew-symmetric, parallel,
torsion were used in \cite{Cleyton1}, as a point of departure to propose a framework somewhat parallel to 
weak holonomy, formulated in terms of $G$-structures (see \cite{Swann1} for a brief account). 
The main result of our paper is a classification result of nearly K\"ahler manifolds, given below.
\begin{teo}
Let $(M^{2n},g,J)$ be a strict and complete, simply connected, nearly K\"ahler manifold. Then $M$ is a 
Riemannian product whose factors belong to one of the following three classes : \\
- $6$ dimensional nearly K\"ahler manifolds; \\
- homogenous nearly K\"ahler spaces of type I,II,III or IV; \\
- twistor spaces over quaternionic K\"ahler manifolds with positive scalar curvature, endowed with the canonical nearly K\"ahler 
metric
\end{teo}
Of course, the previous decomposition coincides with the deRham decomposition of the manifold $(M,g)$. 
The homogenous spaces of type I and II and described in sections 2 and 3 and those of types III and IV in sections 5 and 6. 
The latter are generalized "twistor spaces ", that is Riemannian submersions whose (totally geodesic) fibers are compact 
and simply connected Hermitian symmetric spaces and whose base spaces are compact and simply connected 
symmetric spaces.  \par
Hence, weak holonomy $U(n)$ does not produce any new geometric structure, except possibly in real dimension 
$6$. But even in this case, the only known examples are homogenous. \par
The proof of our theorem goes as follows. In view of the results of \cite{Cleyton1} we are mainly concerned with 
the case when the holonomy of the canonical Hermitian connection is reducible. Then, our point of departure is to 
show that one can suppose that the holonomy representation is of special algebraic type, that is it has strong 
algebraic properties related to the torsion tensor of the Hermitian connection. This will be done in section 3. In the case 
of special algebraic torsion, we show in section 4 that the manifold $M$ is the total space of a fibration whose (totally 
geodesic ) fibers are compact and simply connected Hermitian symmetric spaces. It turns out that for such 
fibrations we can prove an analogue of the DeRham decomposition theorem, in the sense the one can always suppose 
that the fiber and the base space are irreducible, in the usual Riemannian sense. After develloping some facts 
related to the Ricci tensor in section 5 we finally prove the classification result in section 6. It reposes on the 
observation that the torsion tensor points out in a nice way to the Riemannian holonomy group of the base 
manifold, hence allowing the use of the Berger classification theorem.
\par
As a corollary of theorem 1.1 we show how one can use the classification of nearly K\"ahler manifolds in order to study geometric structures on K\"ahler 
manifolds. To cite some related results, recall that it was known for a long time that that a K\"ahler submersion has integrable and totally geodesic horizontal distribution, hence 
zero curvature \cite{Wat1, Johnson}. The same 
conclusion holds for nearly K\"ahler and almost K\"ahler submersions \cite{Falcitelli}. In twistor theory, one searches 
for a non-holomorphic fibration whose total space and fibers are complex manifolds. A reasonable way to see 
these objets arise is from Riemannian foliations. Unfortunately, the following result shows that it cannot happen in the 
presence of a K\"ahler metric, except for a few special cases. 
\begin{teo}
Let (M,g,J) be a simply connected K\"ahler manifold supporting a foliation ${\cal{F}}$ such that $(M,g,{\cal{F}})$ becomes a totally geodesic 
Riemannian foliation. If the leaves of ${\cal{F}}$ are complex manifolds then $(M,g)$ is a Riemannian product  whose factors belong to one the 
following three classes : \\
- twistor spaces over quaternionic K\"ahler manifolds of positive scalar curvature ; \\
- compact homogenous K\"ahler manifolds corresponding to nearly K\"ahler manifolds of types III and IV; \\
- K\"ahler manifolds. \\ 
Furthermore, the foliation ${\cal{F}}$ is obtained in a canonical way from the above decomposition. 
\end{teo}
This is not a very surprising result, as geometric structures on K\"ahler manifolds are rather rare. As an example 
which parallels theorem 1.2 we cite the result that a complex contact manifold with a K\"ahler-Einstein metric 
is the twistor space of a quaternionic K\"ahler manifold of positive sectional curvature \cite{Lebrun1}. The main 
ingredient of the proof of the theorem 1.2, given in section 7, is the observation that K\"ahler manifolds 
admitting a Riemannian foliation as above have a canonical nearly K\"ahler metric. Next we use the results 
leading to the proof of theorem 1.1. \par
Finally, in section 8 we study homogenenous nearly K\"ahler manifolds using theorem 1.1. A classification result involving the previously defined spaces 
of type $I, II, III$ and $IV$ is obtained. As a consequence, we recover by geometric means the 
classification of naturally reductive Riemannian $3$-symmetric spaces, which was 
obtained in his general form in \cite{Gray4} using Lie theory. We are also able to reduce 
a conjecture of Wolf and Gray relating to nearly K\"ahler homogenenous spaces to 
its $6$-dimensional version. \\
{\bf{Acknowledgments.}} The author wishes to thank Andrew Swann for usefull discussions during the preparation of this work and also 
Alain Valette for explanations on representation theory.
$\\$
\section{Preliminaries}
A nearly K\"ahler manifold is an almost Hermitian manifold $(M^{2n},g,J)$ such 
that 
$$ (\nabla_X J)X=0$$
for every vector field $X$ on $M$ (here $\nabla$ denotes the Levi-Civita 
connection associated to metric $g$). It is called {\it{strict}} if $\nabla_XJ \neq 0$ whenever 
$X \in TM, X \neq 0$. \par
Recall that the tensor $\nabla J$ has a number of important 
algebraic properties that can be summarized as follows : the tensors $A$ and $B$ defined for 
$X,Y,Z$ in $TM$ by $A(X,Y,Z)=<(\nabla_XJ)Y,Z>$ and $B(X,Y,Z)=<(\nabla_XJ)Y,JZ>$ are skew-symmetric and 
have type $(0,3)+(3,0)$ as real $3$-forms.
Another object of particular importance is the canonical Hermitian connection defined by 
$$ \overline{\nabla}_XY=\nabla_XY+\frac{1}{2} (\nabla_XJ)JY.$$
It is easy to see that $\overline{\nabla}$ is the unique Hermitian connection on $M$ with totally skew-symmetric 
torsion (see for example \cite{Friedrich2}). Note that the torsion of $\overline{\nabla}$ 
given by $T(X,Y)=(\nabla_XJ)JY$ vanishes iff $(M,g,J)$ is a K\"ahler manifold. \par 
Let us define $\Phi \in \Lambda^3(M)$ by $\Phi(X,Y,Z)=<(\nabla_XJ)Y,Z>$ whenever $X,Y,Z$ 
are in $TM$. A fact of crucial importance for us will be that the form $\Phi$ is parallel with respect 
to the canonical Hermitian connection, that is : 
\begin{nr} \hfill 
$ \overline{\nabla} \Phi=0. \hfill $
\end{nr}
This can be deduced from the $\overline{\nabla}$ parallelism of the tensor $\nabla J$ which 
was proven in $\cite{Moro1}$. This particularly implies that the holonomy group of $\overline{\nabla}$ is contained at each point in the 
$U(n)$-stabiliser of the three form $\Phi$. 
The following decomposition result shows that one can restrict attention to strict nearly K\"ahler manifolds.
\begin{pro} \cite{Kiri, Nagy1}
Let $(M,g,J)$ be a complete and simply connected nearly K\"ahler manifold. Then $M$ is a Riemannian product $M_1 \times M_2$ where 
$M_1$ and $M_2$ are K\"ahler respectively strict nearly K\"ahler manifolds. 
\end{pro}
\par We discuss now briefly 
some properties of the Ricci tensor of nearly K\"ahler manifolds. Recall that it was 
shown in \cite{Nagy1} that if $(M,g,J)$ is strict then it has positive Ricci curvature, hence $(M,g,J)$ is compact as soon as it is complete. Moreover, under these conditions 
the fundamental group  is finite hence the requirement of simple connectivity (to be done in most cases in this work) is not too restrictive. Now, the 
Ricci-$\star$ curvature is defined by : 
$$ Ric^{\star}(X,Y)=\frac{1}{2}\sum \limits_{i=1}^{2n}R(X,JY,e_i,Je_i)$$
where $R$ is the curvature tensor of $(M,g)$ and $\{e_1,\ldots, e_{2n}\}$ a local frame 
field. The difference of the Ricci and Ricci-$\star$ curvature tensors, to be 
denoted by $r$, is given by the formula (see 
\cite{Gray1}) : 
$$ <rX,Y>=\sum \limits_{i=1}^{2n} <(\nabla_{e_i}J)X, (\nabla_{e_i}J)Y>.$$
Then $r$ is symmetric, positive and commutes with $J$; furthermore we have 
$\overline{\nabla}r=0$ (see \cite{Nagy1}). It is well known  \cite{Gray1} that if the tensor $r$ has 
a single eigenvalue then $(M,g)$ is an Einstein manifold and furthermore the first 
Chern class of $(M,J)$ vanishes. \par 
From now on we will work with a strict and complete nearly K\"ahler manifold $(M^{2n},g,J)$.  
The torsion tensor of the canonical Hermitian connection can be used to define a vector 
product on $TM$ by setting 
$$ X \bullet Y=T(X,Y)$$
whenever $X,Y$ belong to $TM$. The algebraic properties of the torsion tensor $T$ 
translate into analogous ones for the previous defined vector product. 
If $V,W$ are distributions of $TM$ we will usually denote by 
$V \bullet W$ the linear span of $\{ x \bullet y : x \in V, y \in W\}$. \par
Let us recall that the first Chern class of $(M,J)$ is represented by the closed $2$-form $\gamma_1$ 
defined by $8 \pi \gamma_1(X,Y)=<CX,JY>$ for all $X,Y$ in $TM$. Here $C$ denotes the symmetric 
endomorphism $Ric-5Ric^{\star}$. It has been proven in \cite{Nagy1} that the 
tensor $C$ is $\overline{\nabla}$-parallel. As $\gamma_1$ is a closed form we obtain easily 
\begin{nr} \hfill 
$ -C(X \bullet Y)=X \bullet CY+CX \bullet Y \hfill$
\end{nr}
for all $X$ and $Y$ in $TM$.
\begin{rema}
A basic fact about $C$ is that if $E_i$ and $E_j$ are eigenspaces of $C$ with corresponding eigenvalues 
$\lambda_i$ and $\lambda_j$ then $-(\lambda_i+\lambda_j)$ is an eigenvalue of $C$ if $E_i \bullet E_j 
\neq 0$.
\end{rema}
\begin{pro} Suppose that $M$ is simply connected. If $C$ has a single eigenvalue then 
it splits as a Riemannian product whose factors are strict nearly K\"ahler manifolds 
such that their corresponding tensors $C$ and $r$ have exactly one eigenvalue. \end{pro}
{\bf{Proof}} : \\
For the inverse implication see \cite{Gray1}. Let us suppose now that $C$ has a single eigenvalue. Then 
$C=0$ by (2.2) and more, $Ric=\frac{5}{4}r$. We recall now the expression of the Ricci curvature of $(M,g)$ 
we obtained in \cite{Nagy1}. If $V_i, 1 \le i \le p$ are the eigenspaces of $r$ corresponding to the 
eigenvalue $\lambda_i$ then 
\begin{nr} \hfill 
$ Ric(X,Y)=\frac{\lambda_i}{4}<X,Y>+\frac{1}{\lambda_i} \sum \limits_{s=1}^{p}\lambda_s<r^s(X),Y>\hfill $
\end{nr}
for all $X$and $Y$ in $V_i$, where the tensors $r^s, 1 \le s \le p$ are defined by $<r^s(X),Y>=-Tr_{V_s}(\nabla_XJ)(\nabla_YJ)$. Let us now 
suppose that the eigenvalues of $r$ are ordered as follows : $0 <\lambda_1 < \ldots < \lambda_p$. Then, on $V_p$ 
we must have $\lambda_p=\frac{1}{\lambda_p}\sum \limits_{s=1}^{p}\lambda_s r^s$. As $r=\sum \limits_{s=1}^p r^s$ 
this implies easily that $r^s=0, 1 \le s \le p-1$ on $V_p$. In other words, $V_p \bullet H=0$ where 
$H= \bigoplus \limits_{s=1}^{p-1}V_s$ and it follows that $V_p \bullet V_p \subseteq V_p$ and $H \bullet H 
\subseteq H$. Hence, the distributions $V_p$ and $H$ are $\nabla$-parallel, thus -by  DeRham's decomposition 
theorem- $M$ is a Riemannian product $M_1 \times M_2$. It easy to show that $M_1, M_2$ are strict 
nearly K\"ahler manifolds and $M_1$ has the property that the corresponding tensors $r$ and $C$ have 
exactly one eigenvalue. Now we apply the same procedure to the manifold $M_2$ and we conclude by induction. q.e.d.
\begin{rema}
(i) Note that a result analogous to proposition 2.1 was already proven by different means in \cite{Wat} under the 
assumptions of Riemannian irreducibility and constant scalar curvature. Note that the latter is now known to hold 
for all strict nearly K\"ahler manifolds. By the work in \cite{VanHecke}, proposition 2.1 remains true if we replace 
the vanishing of the Chern form with the vanishing of the Chern class. \\
(ii) Manifolds with $C=0$ and $r$ with a single eigenvalue are automatically Einstein.
\end{rema}
In the rest of this section we will set up some facts concerning the holonomy representation of the canonical Hermitian connection. We begin with the 
following elementary result in representation theory. We also include its proof as we were unable to find it in the litterature.
\begin{lema} Let $(V^{2n},g,J)$ a real vector space equipped with a scalar product $g$ and a compatible complex structure $J$. Let $V^{\mathbb{C}}$ the complex vector space 
obtained from $V$ by defining the complex multiplication to be $i.v=Jv$. Let $\pi$ be a real representation of a group $G$ on $V$ respecting the inner product and the complex 
structure. Let $\pi^{\mathbb{C}}$ be the complex representation of $G$ on $V^{\mathbb{C}}$ induced by $\pi$. If $\pi^{\mathbb{C}}$ is irreducible then the following possibilities may occur : \\
(i) $\pi$ is irreducible; \\
(ii) there exists an irreducible $G$-submodule of $V$, to be denoted by  ${\cal{V}}$, such that $V={\cal{V}} \oplus J{\cal{V}}$, an orthogonal direct sum.
\end{lema}
{\bf{Proof}} : \\
If $\pi$ is irreducible then we are in the case of (i). Suppose now that the representation of $G$ on $V$ is real reducible. Let ${\cal{V}}$ be a proper invariant subspace. 
Then ${\cal{V}} \cap J{\cal{V}}$ is preserved by $J$ and $G$-invariant so it follows that 
${\cal{V}} \cap J{\cal{V}}=0$. The same argument implies that we have the direct sum decomposition $V={\cal{V}} \oplus J{\cal{V}}$. We will 
show now that the representation of $G$ on ${\cal{V}}$ is irreducible. Indeed, if $E \subseteq {\cal{V}}$ is 
non-zero $G$ module then as before $E \cap JE=0$ and $E \oplus JE=V$, hence $E$ and ${\cal{V}}$ have 
the same dimension and thus $E={\cal{V}}$. It remains to show that the sum ${\cal{V}} \oplus J{\cal{V}}$ is orthogonal. Using Riesz's representation theorem we obtain 
the existence of a skew-symmetric endomorphism of ${\cal{V}}_m$, to be denoted by $F$, such that : 
$$ <v,Jw>=<Fv,w>  \ \mbox{for all } \ v,w \ \mbox{in} \ {\cal{V}}.$$
Of course, $F$ must be $G$-invariant. We extend $F$ to a endomorphism $f$ of $V$, by setting $f(Jv)=J(Fv)$ for all $v$ in ${\cal{V}}$. Obviously this is complex linear and $G$-invariant, 
so using the irreducibility of $V^{\mathbb{C}}$  and Schur's lemma we obtain that $f=z1_{V^{\mathbb{C}}}$ for some complex number $z$. The definition of $f$ and the skew-symmetry of $F$ 
imply now the vanishing of $F$. q.e.d.
\par
Before applying  this to our geometric context let us recall the following important result. 
\begin{teo} \cite{Cleyton1}
Let $(M^n,g)$ be a simply connected Riemannian manifold equipped with a metric connexion $\nabla^{\prime}$ whose torsion tensor, denoted by $T$, is 
$\nabla^{\prime}$-parallel and totally skew. If the Lie algebra of the stabiliser of $T$ in $\mathfrak{so}(n)$ acts 
irreducibly on some tangent space, then $M$ is a homogenous space with irreducible isotropy unless $Stab_{\mathfrak{so}(n)}T$ equals 
$G_2$ or $SU(3)$ acting as usually on $\mathbb{R}^7$ and $\mathbb{R}^6$ respectively.
\end{teo}
\begin{rema}
The two exceptions in theorem 2.1 are weak holonomy structures, precisely nearly parallel $G_2$ and six dimensional nearly 
K\"ahler manifolds.The basic ingredient for the proof of the previous result is a structure theorem for the space of algebraic curvature 
tensors of a Lie algebra preserving a three form. The homogenous structure on $M$ is obtained by proving that $\nabla^{\prime}$ is an 
Ambrose-Singer connexion, that is $\nabla^{\prime} R^{\prime}=0$. Also, it would be interesting to have a differential geometric proof of this result.
\end{rema}
This remark motivates the setting, in the nearly K\"ahler case, of the following definition. 
\begin{defi}
A strict nearly K\"ahler manifold $(M^{2n},g,J)$ with $n \ge 4$ is said homogenous of type I if the holonomy representation of the canonical Hermitian connection is 
real irreducible.
\end{defi}
By the previous discussion, in this case $\overline{\nabla}$ is a Ambrose-Singer connexion and moreover, $M$ is isotropy irreducible. We can now use the previous preliminaries 
in order to distinguish nearly K\"ahler manifolds by holonomic means, as follows.
\begin{pro}
Let $(M^{2n},g,J)$ a strict, complete and simply connected nearly K\"ahler manifold. Then the following possibilities may occur : \\
(i) $M$ is homogenous of type I, unless $M$ is of dimension $6$ ; \\
(ii) we have a $\overline{\nabla}$-parallel, orthogonal, decomposition $TM={\cal{V}} \oplus H$ with $H=J{\cal{V}}$ ;  \\
(iii) the holonomy representation of the canonical Hermitian connection is complex reducible. 
\end{pro}
{\bf{Proof}} : \\
We use that the torsion of the canonical Hermitian connection is parallel and totally skew. In particular, this implies 
that if the holonomy representation of $\overline{\nabla}$ is real irreducible then the same is true for representation of the $\mathfrak{so}(2n)$-stabilizer of the 
three form $\Phi$. We conclude by lemma 2.1 and theorem 2.1. q.e.d.
\par
Notice the difference with the classical Riemannian case, where there is no scarcity of holonomy irreducible structures. 
We mention apart the coresponding version of proposition 2.3 in six dimensions. We first set first the following 
\begin{defi}
Let $(M,g,J)$ be a $6$-dimensional nearly K\"ahler manifold. It is called Hermitian irreducible if the holonomy representation 
of the canonical Hermitian connexion is irreducible in the real sense (and hence the Hermitian holonomy equals $SU(3)$).
\end{defi}
\begin{coro}
Let $(M,g,J)$ be a complete and simply connected nearly K\"ahler $6$-fold. The following cases may occur : \\
(i) $(M,g,J)$  is Hermitian irreducible; \\
(ii) we have a $\overline{\nabla}$-parallel, orthogonal, decomposition $TM={\cal{V}} \oplus H$ with $H=J{\cal{V}}$ ;  \\
(iii) (M,g,J) is one of the spaces $P\mathbb{C}^3, \mathbb{F}^3$ equipped with its canonical nearly K\"ahler metric.
\end{coro}
{\bf{Proof}} : \\
Follows immediately from proposition 2.3 and \cite{Moro1} (see also \cite{Nagy1}). q.e.d.
\par
Thus, in order to achieve the classification of strict nearly K\"ahler manifolds 
it remains to treat the cases when the holonomy representation of the canonical Hermitian connection 
is reducible, first in the real sense precised at the point (ii) of the proposition 2.3 and next in the general complex sense. We leave open the problem of 
classifying Hermitian irreducible nearly K\"ahler $6$-folds.

\section{Curvature properties and reducibility}
In this section we will consider a strict nearly K\"ahler manifold $(M^{2n},g,J)$ which is reducible in the real 
Hermitian sense, that is the tangent bundle of $M$ admits a $\overline{\nabla}$-parallel, orthogonal, decomposition : 
$$ TM={\cal{V}} \oplus H. $$
In the first part of this section we will show that it is always possible to suppose that the previous 
decomposition is of special type, that is ${\cal{V}} \bullet {\cal{V}} \subseteq {\cal{V}}$. This will 
be done by a carefull examination of the curvature 
tensor of the canonical Hermitian connection, to be denoted by $\overline{R}$. Some preliminaries are required.
\begin{lema}
Let $X$ be in $H, Y$ in $TM$ and $V, W$ in ${\cal{V}}$. We have : 
$$\overline{R}(X,Y,V,W)=<[\nabla_{V}J, \nabla_{W}J]X,Y>-<(\nabla_XJ)Y, (\nabla_{V}J)W>. $$ 
\end{lema}
{\bf{Proof}} : \\
Let us recall the relation (see \cite{Gray1} page 237) : 
\begin{nr} \hfill
$\begin{array}{lr} 
\overline{R}(X,Y,Z,T)=R(X,Y,Z,T)-\frac{1}{2}<(\nabla_XJ)Y, 
(\nabla_ZJ)T>+ \\
\frac{1}{4} \biggl [ <(\nabla_XJ)Z, (\nabla_YJ)T>-<(\nabla_XJ)T, (\nabla_YJ)Z> \biggr ] 
\end{array} \hfill $
\end{nr}
where $R$ denotes the curvature tensor of the Levi-Civita connection of $g$. It follows by the first Bianchi identity that 
$$ R(X,Y,V_1,V_2)+R(Y,V_1,X,V_2)+R(V_1,X,Y,V_2)=0.$$
Or $\overline{R}(Y,V_1,X,V_2)=0$ so by (3.1) we get 
$$ \begin{array}{lr}
R(Y,V_1,X,V_2)=\frac{1}{2}<(\nabla_YJ)V_1, (\nabla_XJ)V_2> \\
-\frac{1}{4} \biggl [<(\nabla_YJ)X,(\nabla_{V_1}J)V_2>-
<(\nabla_YJ)V_2,(\nabla_{V_1}J)X> \biggr ].
\end{array}$$
In the same manner we obtain an analogous formula for $R(V_1,X,Y,V_2)$ and the result 
follows easily.  q.e.d. \\ \par
\begin{coro}
Let $(M^{2n},g,J)$ be complete, strict and simply connected nearly K\"ahler manifold satisfying the conditions of (ii) in proposition 2.3. Then
$(M,g)$ is a homogenous space.
\end{coro}
{\bf{Proof}} : \\
Let $V_i, 1 \le i \le 4$ be in ${\cal{V}}$. Then using lemma 3.1 we get 
\begin{nr} \hfill
$ \overline{R}(V_1,V_2,JV_3,JV_4)=<[\nabla_{V_1}J,
\nabla_{V_2}J]V_3,V_4>+<(\nabla_{V_1}J)V_2, (\nabla_{V_3}J)V_4>. \hfill $
\end{nr}
Using now that $\overline{R}(V_1,V_2,JV_3,JV_4)=\overline{R}(V_1,V_2,V_3,V_4)=
\overline{R}(JV_1,JV_2,JV_3,JV_4)$ it follows that (3.2) completely determines the curvature $\overline{R}$ on 
$TM$. We recall 
now that $\overline{\nabla} (\nabla J)=0$, a fact proven in \cite{Moro1}, which is equivalent to 
\begin{nr} \hfill 
$ \overline{\nabla}_U((\nabla_XJ)Y)=(\nabla_XJ)\overline{\nabla}_UX+(\nabla_{\overline{\nabla}_UX}J)Y \hfill $
\end{nr}
whenever $U,X,Y$ are in $TM$. Thus, derivating formulas of type (3.2) we get easily that $\overline{\nabla} \ 
\overline{R}=0$, hence $\overline{\nabla}$ is an Ambrose-Singer connection. As $M$ is simply connected, we have 
that $(M,g,J)$ is a homogenous space (see \cite{Tri}).  q.e.d.
\\ \par
It is convenient to distinguish now a second class of homogenous nearly K\"ahler manifolds as 
follows.
\begin{defi}
A simply connected, strict, nearly K\"ahler manifold is homogenous of type II iff the holonomy representation of $\overline{\nabla}$
on $T_m^{\mathbb{C}}(M)$ is complex irreducible and satisfies the conditions of (ii) in proposition 2.3.
\end{defi}
\begin{rema}
(i) In section 8 we will study in more detail the class of nearly K\"ahler manifolds having the property that $\overline{\nabla}$ is an 
Ambrose-Singer connection. \\
(ii) A typical example of nearly K\"ahler manifold of type II is provided by products $G \times G$  where $G$ is a compact Lie group 
without center (see \cite{Gray4}). However, it is not straightforward to see that these examples exhaust our class of spaces. This problem 
is studied in \cite{Nagy3}. \\
\end{rema}
Hence, when the holonomy representation of the canonical Hermitian connexion is complex irreducible, $(M,g,J)$ is homogeneous of type I or II. 
Therefore, it remains 
to understand the complex reducible case. In this situation, using parallel transport 
we get a $\overline{\nabla}$-parallel decomposition $TM={\cal{V}} \oplus H$ which 
is {\it{stable}} by $J$. For this reason, all the distributions we will consider from now on will be preserved by $J$.\par
Now, let us state and prove the following important consequence of lemma 3.1.
\begin{coro} Let $X,Y$  be in $H$ and $V,W$ be in ${\cal{V}}$ respectively. Then : \\
(i) $(\nabla_XJ)(\nabla_VJ)W=0$; \\
(ii) $(\nabla_XJ)(\nabla_YJ)Z$ belongs to $H$ whenever $Z$ belongs to $H$; \\
(iii) $(\nabla_VJ)(\nabla_WJ)X$ belongs to $H$; \\
(iv) $(\nabla_XJ)(\nabla_YJ)V$ belongs to ${\cal{V}}$.
\end{coro}
{\bf{Proof}} : \\
(i) We use that $\overline{R}(JX,JZ,V,W)=\overline{R}(X,Z,V,W)$ for all $Z$ in $TM$, lemma 3.1 and the algebraic properties of 
the tensor $\nabla J$ (see Section 2).\\
(ii) Reversing the roles of ${\cal{V}}$ and $H$ we obtain that $(\nabla_VJ)(\nabla_XJ)Y=0$ for all $X,Y$ in $H$ and 
$V $ in ${\cal{V}}$. Taking the scalar product with $Z$ in $H$ gives now the result. \\
(iii) We choose $Y=V_1$ in lemma 3.1, where $V_1$ belongs to ${\cal{V}}$. Since $\overline{R}(X,V_1,V,W)=0$ we get 
by (i) that $[\nabla_VJ,\nabla_WJ]X$ belongs to $H$. But $[\nabla_VJ, \nabla_{JW}J]JX=-\{\nabla_VJ, \nabla_WJ\}X$ 
is equally in $H$ and the result follows. \\
(iv) It suffices to interchange the roles of the distributions ${\cal{V}}$ and $H$. q.e.d.
\\ \par
We can now show the following result which mainly asserts that the tangent bundle of a reducible 
nearly K\"ahler manifold contains an integrable distribution. 
\begin{pro}
Let $(M^{2n},g,J)$ be a complete, strict nearly K\"ahler manifold and let us suppose that $M$ is complex reducible and simply connected. 
Then $M$ splits as a Riemannian product $Z \times N$ 
where $Z$ and $N$ are complete and strict nearly K\"ahler manifolds and the tangent bundle of $N$ 
contains a proper $\overline{\nabla}$-parallel distribution with ${\cal{V}}$ with ${\cal{V}} \bullet {\cal{V}}=0 $.
\end{pro}
{\bf{Proof}} : \\
We will show first that one can find a $\overline{\nabla}$-parallel distribution ${\cal{V}}$ in $TM$ such that 
${\cal{V}} \bullet {\cal{V}} \subseteq {\cal{V}}$. Indeed, 
the reducibility of $M$ implies the existence of a $\overline{\nabla}$-parallel decomposition $TM=
E \oplus F$. Let $F_0$ be the distribution generated by elements of the form $(v \bullet w)_F$ where 
$v,w$ belong to $E$ and the subscript denotes orthogonal projection on $F$. Using corollary 3.2, (i) we get 
$$ (\nabla_XJ)(\nabla_YJ) (v \bullet w)_E +(\nabla_XJ)(\nabla_YJ) (v \bullet w)_F=0$$
whenever $X,Y$ are in $F$. But the first term of this sum belongs to $E$ by corollary 3.2, (iv) while the 
second is in $F$ by corollary 3.2, (ii). It follows that each term of the previous sum vanishes and this implies 
easily that $F \bullet F_0=0$. As $F_0$ is contained in $F$ we get that $F_0 \bullet F_0=0$. Obviously, $F_0$ 
is $\overline{\nabla}$-parallel (one uses (3.3)) hence we may take ${\cal{V}}=F_0$ when $F_0 \neq 0$. If $F_0$ then 
$E \bullet E \subseteq E$ and we set ${\cal{V}}=E$. \par
Therefore, let us consider a $\overline{\nabla}$-parallel decomposition $TM={\cal{V}} \oplus H$, such that 
${\cal{V}}\bullet  {\cal{V}}\subseteq  {\cal{V}} $.
Let $r_1 : V \to V$ be the tensor defined by $<r_1V,W>=-Tr_{{\cal{V}}} (\nabla_VJ)(\nabla_WJ)$ for all 
$V,W$ in ${\cal{V}}$. Then $r_1$ is parallel with respect to the unitary connection induced by $\overline{\nabla}$ 
in the bundle ${\cal{V}}$ thus we have the decomposition ${\cal{V}}={\cal{V}}_0 \oplus {\cal{V}}_1$ with 
${\cal{V}}_0=Ker(r_1)$ and ${\cal{V}}_1$ the 
orthogonal complement of ${\cal{V}}_0$ in ${\cal{V}}$.  Since the corollary 3.2, (i) implies 
that $(\nabla_XJ)r_1W=0$ for all $X$ in $H$ and $W$ in ${\cal{V}}$ it follows that $H \bullet {\cal{V}}_1=0$. Thus, we have 
the $\overline{\nabla}$ parallel decomposition 
$$ TM={\cal{V}}_1 \oplus H_1$$
where $H_1={\cal{V}}_0 \oplus H$. Since the very definition of ${\cal{V}}_0$ implies that 
${\cal{V}} \bullet {\cal{V}}_0=0$ it is now easy to establish that 
$H_1 \bullet H_1 \subseteq H_1, {\cal{V}}_1 \bullet {\cal{V}}_1 \subseteq {\cal{V}}_1$ and 
${\cal{V}}_1 \bullet H_1=0$. It follows that the distributions ${\cal{V}}_1$ and 
$H_1$ are in fact $\nabla$-parallel and we conclude by the decomposition theorem of de Rham.  q.e.d.
\begin{rema}
In the case when $C \neq 0$ (see section 2 for definitions) proposition 3.1 can be given a simple algebraic proof.
Let $\lambda_i, 1 \le i \le p$ be the (pairwise distinct) eigenvalues of $C$ ordered such that $ \vert \lambda_1 \vert \le 
\ldots \le \vert \lambda_p \vert$. If $E_i,  1 \le i \le p$ are the corresponding eigenspaces using the remark 2.1 it follows 
that $E_p \bullet E_p=0$. \\
Nevertheless, in the case when $C=0$ (which implies that $M$ is Einstein) the 
geometric arguments of the proof of the proposition 3.1 are needed.
\end{rema}
We will show now that the class of nearly K\"ahler manifolds appearing in proposition 3.1 leads to 
another class introduced by the following definition.
\begin{defi} Let $(M^{2n},g,J)$ be a strict and complete nearly K\"ahler manifold. It is said to have 
special algebraic torsion if there exists a $\overline{\nabla}$-parallel decomposition $TM={\cal{V}} \oplus H$ with 
${\cal{V}} \bullet {\cal{V}}=0$ and $H \bullet H={\cal{V}}$.
\end{defi}
\begin{rema}
A strict nearly K\"ahler manifold (M,g,J) admitting a $\overline{\nabla}$-parallel decomposition $TM={\cal{V}} \oplus H$ where 
${\cal{V}} \bullet {\cal{V}}=0$ et $H \bullet H \subseteq {\cal{V}}$ is of special algebraic type.  To see this, let us set 
${\cal{V}}_0=H \bullet H$ and let ${\cal{V}}_1$ be the orthogonal complement of 
${\cal{V}}_0$ in ${\cal{V}}$. By the definition of ${\cal{V}}_1$ we must have ${\cal{V}}_1 \bullet H=0$ hence 
${\cal{V}}_1 \bullet TM=0$ as 
${\cal{V}} \bullet {\cal{V}}=0$. The fact that $M$ is strict implies now that ${\cal{V}}_1=0$, thus $H \bullet H={\cal{V}}$. Using a similar argument 
one can also prove that ${\cal{V}} \bullet H=H$.
\end{rema}
\begin{pro}
Let $(M^{2n},g,J)$ be a nearly K\"ahler manifold equipped with a 
$\overline{\nabla}$-parallel distribution ${\cal{V}}$ satisfaying 
${\cal{V}} \bullet {\cal{V}}=0$. If $M$  is simply connected, it splits as a Riemannian 
product $Z \times N$ where both factors are strict nearly K\"ahler manifolds and $N$ has 
special algebraic torsion.
\end{pro}
{\bf{Proof}} : \\
Set $H_0={\cal{V}} \bullet H$ which is obviously contained in $H$. If $H_1$ is the orthogonal complement of 
$H_0$ in $H$ we get a $J$ invariant, $\overline{\nabla}$-parallel decomposition $H=H_0 \oplus H_1$. Then 
${\cal{V}} \bullet H_0 \subseteq H_0$ and ${\cal{V}} \bullet H_1 \subseteq H_0$ by the definition of $H_0$. But then 
${\cal{V}} \bullet H_1$ is orthogonal to $H_0$ hence ${\cal{V}} \bullet H_1=0$. Using corollary 3.2, (iii) we obtain 
that $H \bullet H_0 \subseteq {\cal{V}}$ and since ${\cal{V}} \bullet H_1=0$ we must have $H_0 \bullet H_1=0$. So 
$H_0 \bullet H_0 \subseteq {\cal{V}}$ and it is now easy to see that $H_1 \bullet H_1 \subseteq H_1$. Hence the 
decomposition $TM=E_1 \oplus E_2$  is $\overline{\nabla}$-parallel with $E_i \bullet E_i \subseteq E_i, i=1,2, 
E_1 \bullet E_2$  and we conclude using the de Rham decomposition theorem and the previous remark.  q.e.d.
\\ \par
The results of this section lead easily to the following theorem which roughly classifies strict nearly K\"ahler manifolds 
up to those having special algebraic torsion. 
\begin{teo}
Let $(M^{2n},g,J)$ be a complete, simply connected and strict nearly K\"ahler manifold. Then $M$
decomposes as Riemannian product whose factors belong to one of the following three classes of 
strict nearly K\"ahler manifolds : \\
-homogenous spaces of type I and II; \\
-Hermitian irreducible $6$-dimensional manifolds; \\
-manifolds with special algebraic torsion.
\end{teo}
\section{Special algebraic torsion}
Let $(M,g,J)$ be nearly K\"ahler with special algebraic torsion. Thus, we have a $\overline{\nabla}$-parallel 
decomposition $TM={\cal{V}} \oplus H$ with ${\cal{V}} \bullet {\cal{V}}=0$ and $H \bullet H = {\cal{V}}$. Our 
first observation is that ${\cal{V}}$ being 
$\overline{\nabla}$-parallel with ${\cal{V}} \bullet {\cal{V}}=0$ it is integrable. We would like to establish that 
${\cal{V}}$ induces on $M$ the structure of a fibration, with smooth base space. We first prove the following result, showing that the torsion tensor of 
the canonical Hermitian connection completely determines 
the curvature in the ${\cal{V}}$-direction.
\begin{pro}
Let $X,Y$ be in $H$ and $V_1, V_2,V_3$ in $V$. Then :
 \begin{nr} \hfill
$ \overline{R}((\nabla_XJ)JY, V_1,V_2,V_3)=
<JY, [\nabla_{V_1}J, [\nabla_{V_2}J,\nabla_{V_3}J]]X>. \hfill $
\end{nr}
Moreover, we have $\overline{\nabla}_U \overline{R}(V_1,V_2,V_3,V_4)=0$.
\end{pro}
{\bf{Proof}} : \\
Using the second Bianchi identity for the Hermitian connection $\overline{\nabla}$ (see \cite{KN}) we obtain  
$$ \begin{array}{lr}
(\overline{\nabla}_X \overline{R})(Y,V_1,V_2,V_3)+(\overline{\nabla}_Y\overline{R})(V_1,X,V_2,V_3)+
(\overline{\nabla}_{V_1}\overline{R})(X,Y,V_2,V_3)+\vspace{2mm}\\
\overline{R}((\nabla_XJ)JY, V_1,V_2,V_3)+\overline{R}((\nabla_YJ)JV_1, X, V_2,V_3)+
\overline{R}((\nabla_{V_1}J)JX, Y, V_2,V_3)=0.
\end{array}$$
Since the distributions ${\cal{V}}$ and $H$ are $\overline{\nabla}$-parallel we get 
$(\overline{\nabla}_X \overline{R})(Y,V_1,V_2,V_3)=\\ (\overline{\nabla}_Y\overline{R})(V_1,X,V_2,V_3)=0$ and 
using lemma 3.1 and (3.3) it can be seen that the term
$(\overline{\nabla}_{V_1}\overline{R})(X,Y,V_2,V_3)$ also vanishes. Now, we compute the last two curvature 
terms in the Bianchi identity using lemma 3.1 and the corollary 3.1, (i) and the result follows by calculus. Let us prove now 
the second assertion. We first remark that formula (4.1) can be given the form : 
\begin{nr} \hfill 
$ (\nabla_XJ) \overline{R}(V_1,V_2)V_3=[\nabla_{V_3}J, [\nabla_{V_1}J, \nabla_{V_2}]]X\hfill $
\end{nr}
Derivating this in the direction $U$ in $TM$, we get by using (3.3) that 
$$<(\nabla_XJ) \biggl [ \overline{\nabla}_U \overline{R}(V_1,V_2,V_3) \biggr ],Y>=0$$
and we conclude using the special algebraic torsion assumptions, and the strictness of $M$. q.e.d.
\\ \par
As $H \bullet H={\cal{V}}$ it follows that the restriction of the tensor $\overline{R}$ to the distribution ${\cal{V}}$ is completely determined by the formula (4.1). We can now prove 
the following.
\begin{pro} 
There exists a Riemannian manifold $(N,h)$ and a Riemannian submersion with 
totally geodesic fibers $ \pi : M \to N$ whose vertical distribution equals ${\cal{V}}$. With respect to 
the induced metric and almost complex structures, each fiber is a simply connected, compact, Hermitian symmetric space. 
\end{pro}
{\bf{Proof}} : \\
Consider the foliation induced by the distribution ${\cal{V}}$. One can easily show that ${\cal{V}}$ is totally geodesic hence 
we obtain by the second part of the proposition 4.1 that each leaf is a Hermitian symmetric space.
We look now at the Ricci curvature of the leaves. Let $\widehat{Ric} : {\cal{V}} \to {\cal{V}}$ be the Ricci curvature 
in the vertical direction. We consider $V$ be in ${\cal{V}}$ with $V \neq 0$ and let 
$\{e_i\}$ be a local orthonormal basis in ${\cal{V}}$. Taking $V_1=V_3=e_i, V_2=V$ in (4.2), we obtain 
after taking the trace over $H$ that : 
$$<\widehat{Ric}(V),r(V)>=2 \sum \limits_{e_i \in {\cal{V}}}^{} \Vert (\nabla_VJ)(\nabla_{e_i}J)\Vert^2 >0$$
(here, we considered that the norm of the linear operator $A$ on $TM$ is given by $\Vert A \Vert^2=Tr(AA^{\star})$).
Now the tensor $r$ preserves ${\cal{V}}$ and its restriction to ${\cal{V}}$ is $\nabla$-parallel with strictly 
positive eigenvalues. It follows then easily that each leaf has positive Ricci curvature 
and we use a result of \cite{Kob1} to obtain that the leaves are compact and simply connected. The simple 
connectivity of the leaves implies that the foliation has trivial leaf holonomy, so using the criterion from 
\cite{Sharpe}, page 90, we get a smooth fibration $\pi : M \to N$. The metric part of our statement is a standard one, 
which is left to the reader. q.e.d.
 \\ \par 
We will refer to the fibration $\pi : M \to N$ as the {\it{canonical}} fibration of the nearly K\"ahler manifold 
with special algebraic torsion $M$. Its generic fiber will be denoted by $F$. Before going further, let us examine a few geometric properties 
related to special algebraic torsion.
\begin{pro}
Let $(M,g,J)$ be a nearly K\"ahler with special algebraic torsion and corresponding decomposition 
$TM= {\cal{V}} \oplus H$. Define a Riemannian metric $\overline{g}$ on $M$  by  $\overline{g}
=2g_{{\cal{V}}} \oplus g_H$ where $g_{{\cal{V}}}$ and $g_H$ are the restrictions of the metric $g$ to the 
distributions ${\cal{V}}$ and $H$ respectively. We also define a new almost complex structure $\overline{J}$ 
by setting $\overline{J}_{{\cal{V}}}=-J$ and $\overline{J}_{\vert H}=J$. Then $(M,\overline{g}, \overline{J})$ is a K\"ahler manifold of positive 
Ricci curvature, and hence simply connected.
\end{pro}
{\bf{Proof}} : \\
The proof of the fact that $(M, \overline{g},\overline{J})$ 
is K\"ahler is essentially the same with the one given in the $6$-dimensional case in \cite{Moro1}. 
Let us denote by $\widetilde{\nabla}$ the Levi-Civita connection of the metric $\overline{g}$. In the standard 
way we get 
$$\begin{array}{lr}
\widetilde{\nabla}_XY=\nabla_XY, \ \widetilde{\nabla}_XV=\overline{\nabla}_XV-\frac{1}{2}(\nabla_XJ)JV \\
\widetilde{\nabla}_VX=\overline{\nabla}_VX, \ \widetilde{\nabla}_VW=\nabla_VW
\end{array}$$
whenever $X,Y$ are in $H$ and $V,W$ in ${\cal{V}}$. After a routine computation we obtain 
\begin{nr} \hfill 
$ \begin{array}{lr}
Ric_{\widetilde{g}}=Ric_g-\frac{r}{4} \ \mbox{on} \ H \\
Ric_{\widetilde{g}}=Ric_g+\frac{3r}{4} \ \mbox{on} \ {\cal{V}}
\end{array} \hfill $ 
\end{nr}
and of course $Ric_{\overline{g}}(V,X)=0$ for $V,X$ belonging to ${\cal{V}}$ and $H$ respectively. Using (2.3) we 
see that the Ricci tensor of $\widetilde{g}$ is strictly positive and by a theorem of Kobayashi \cite{Kob1} $M$ is 
simply connected. q.e.d.

We are now going to prove that the standard de Rham decomposition of the fiber induces a 
splitting at the level of the canonical fibration.
\begin{teo}Let $(M, g,J)$ be a nearly K\"ahler manifold with special algebraic torsion. Then $M$ is a 
Riemannian product 
$$M_1 \times M_2 \times \ldots \times M_q$$ 
where each $M_i, 1 \le i \le q$, is a nearly K\"ahler manifold with special algebraic torsion and such 
that the canonical fibration has irreducible fiber.
\end{teo}
{\bf{Proof}} : \\
If the fiber $F$ is irreducible there is nothing to prove. Otherwise it decomposes as $F_1 \times F_2$ and 
using parallel transport we obtain a $\nabla$-parallel decomposition ${\cal{V}}={\cal{V}}^1 \oplus {\cal{V}}^2$.
Hence $R((\nabla_XJ)Y,V,V_1,V_2)=0$ for all $X,Y$ in $H, V$ in ${\cal{V}}$ and $V_i $ in ${\cal{V}}^i, i=1,2$. Using 
(4.1) we obtain 
\begin{nr} \hfill 
$ [\nabla_VJ, [\nabla_{V_1}J, \nabla_{V_2}J]]=0.\hfill $
\end{nr}
Taking  $V=V_1$ we get : 
$$ (\nabla_{V_1}J)^2(\nabla_{V_2}J)+(\nabla_{V_2}J)(\nabla_{V_1}J)^2=2(\nabla_{V_1}J)
(\nabla_{V_2}J)(\nabla_{V_1}J).$$
We change $V_2$ in $JV_2$ in the previous relation and we use that $(\nabla_{JV_2}J)=(\nabla_{V_2}J)J$. We get 
$$ (\nabla_{V_1}J)^2(\nabla_{V_2}J)+(\nabla_{V_2}J)(\nabla_{V_1}J)^2=-2(\nabla_{V_1}J)
\nabla_{V_2}J)(\nabla_{V_1}J)$$
hence we must have 
$$ 
(\nabla_{V_1}J)^2(\nabla_{V_2}J)+(\nabla_{V_2}J)(\nabla_{V_1}J)^2=(\nabla_{V_1}J)
(\nabla_{V_2}J)(\nabla_{V_1}J)=0.$$
This implies that $(\nabla_{V_1}J)^3(\nabla_{V_2})=0$ and since $-(\nabla_{V_1}J)^2$ is a nonnegative 
operator we have that $(\nabla_{V_1}J)^2(\nabla_{V_2}J)=0$. We equally have 
$(\nabla_{V_2}J)(\nabla_{V_1}J)^2=0$. It follows that 
$$ \Vert ((\nabla_{V_2}J)(\nabla_{V_1}J))X\Vert^2=\Vert ((\nabla_{V_1}J)(\nabla_{V_2}J))X\Vert^2=0$$
for all $X$ in $H$, that is 
\begin{nr} \hfill 
$(\nabla_{V_2}J)(\nabla_{V_1}J))=(\nabla_{V_1}J)(\nabla_{V_2}J))=0 \hfill $
\end{nr}
whenever $V_1$ and $V_2$ belong to ${\cal{V}}_1$  and ${\cal{V}}_2$ respectively. We define now 
the subspaces $H^1$ and $H^2$ of $H$ by 
$$ H^i={\cal{V}}^i \bullet H, \ i=1,2. $$
Using (4.4) and the fact that ${\cal{V}} \bullet H=H$ (see remark 3.3) we easily obtain that $H=H^1 \oplus H^2$. After standard manipulations it follows that 
${\cal{V}}^1 \bullet H^2={\cal{V}}^2 \bullet H^1=H^1 \bullet H^2 =0$ and further that 
${\cal{V}}^i \bullet H^i =H^i, H^i \bullet H^i ={\cal{V}}^i$ for $i=1,2$.
We are going to show now that the four distributions ${\cal{V}}^1, {\cal{V}}^2, H^1, H^2$ are 
$\overline{\nabla}$-parallel. Using that  ${\cal{V}}^1 \bullet H^1=H^1$ and (3.3) we establish that $H^1$ is 
$\overline{\nabla}$-parallel. For the parallelism of ${\cal{V}}^1$ we use that $H^1 \bullet H^1={\cal{V}}^1$ and 
the remaining cases are analogous. Now, the distributions $E_i={\cal{V}}^i \oplus H^i , i=1,2$ are 
$\overline{\nabla}$-parallel. Since $E_i \bullet E_i \subseteq E_i, i=1,2$ and $E_1 \bullet E_2=0$ they are in fact $\nabla$-parallel and we 
conclude by the de Rham's theorem, as $M$ is simply connected. q.e.d.
\\ \par 
We can therefore admit without loss of generality that the fiber $F$ of the canonical fibration is an irreducible Riemannian manifold. This assumption is to be made in the 
rest of the present section. 
\par
Let us define a tensor ${\cal{R}} : H \times H \to End_{-}(H)$ by 
$$ {\cal{R}}(X,Y)=\overline{R}(X,Y)+\nabla_{(\nabla_XJ)Y}$$
for all $X,Y$ in $H$. Now, if $X^{\star}$ denotes the horizontal lift of the vector field $X$ of $N$, using standard computations from \cite{Besse1}, page 241, we get 
\begin{nr} \hfill 
$ {\cal{R}}(X^{\star},Y^{\star})Z^{\star}=(R^h(X,Y)Z)^{\star}\hfill $
\end{nr}
whenever $X,Y,Z$ are vector fields on $N$ (here $R^h$ denotes the curvature tensor of $N$). 
\begin{lema}
We have 
$$ \sum \limits_{e_i \in H}^{}{\cal{R}}(e_i, (\nabla_VJ)e_i)=\nabla_{r(V)}J$$
for all $V$ in ${\cal{V}}$ and all orthonormal basis $\{ e_i \}$ of $H$. 
\end{lema}
{\bf{Proof}} : \\
We compute the given sum using the base $\{Je_i \}$ and we take into account that $\overline{R}(JX,JY)=\overline{R}(X,Y)$ whenever 
$X,Y$ belong to $H$.  q.e.d.
 \\ \par
\begin{teo}
The Riemannian manifold $(N,h)$ is irreducible.
\end{teo}
{\bf{}Proof} : \\
If $(M,h)$ wasn't irreducible its tangent bundle admitted a decomposition $TN=E_1 \oplus E_2$ parallel with respect of the Levi-Civita 
connection of $h$, to be denoted by $\nabla^h$. We obviously have $H=H_1 \oplus H_2$ where$$ H_i =\{ v \in H : (d \pi )v \in E_i\}$$
for $i=1,2$. Now, the operator $R^h(X,Y)$ preserves $E_1$ and $E_2$ for all $X,Y$ in $TN$. It follows by (4.5) that ${\cal{R}}(X,Y)$ preserves $H_1$ and $H_2$ 
whenever $X$ and $Y$ are in $H$. Hence lemma 4.1 easily implies that ${\cal{V}} \bullet H_i \subseteq H_i$ for $i=1,2$. It follows that $H_1 \bullet H_2$ and 
${\cal{V}}$ are orthogonal but the fact that $H_1 \bullet H_2 \subseteq {\cal{V}}$ gives that $H_1 \bullet H_2=0$. We get a decomposition 
$$ TM={\cal{V}}^{\prime} \oplus H_2 $$
where ${\cal{V}}^{\prime}={\cal{V}}^{\prime} \oplus H_1$. This decomposition satisfies ${\cal{V}}^{\prime} \bullet {\cal{V}}^{\prime} \subseteq {\cal{V}}^{\prime}$ and 
$H_2 \bullet H_2 \subseteq {\cal{V}}^{\prime}$. Let us show now that it is also $\overline{\nabla}$-parallel. \par
As $\overline{\nabla}_{X^{\star}}Y^{\star}=(\nabla_XY)^{\star}$ for all $X,Y$ in $TN$ it follows that $\overline{\nabla}_XY$ belongs to $H_i$ whenever $X$ is in $H$ and 
$Y$ is in $H_i, i=1,2$. It remains to see that $\overline{\nabla}_VX$ belongs to $H_i$ whenever $X$ is in $H_i$, $V$ is in ${\cal{V}}$ and $i=1,2$. Indeed if $X$ is in $E_1$ we know 
that $[V, X^{\star}]$ is in ${\cal{V}}$ so by the parallelism of $H$ is follows that $\overline{\nabla}_V X^{\star}=-(\nabla_VJ)JX^{\star}$ which belongs to $H_1$ as we already shown 
that ${\cal{V}} \bullet H_1 \subseteq H_1$. Thus, $H_1$ is $\overline{\nabla}$-parallel and of course the same holds for $H_2$. \par
Hence we may apply corollary 3.2, (i) for the distribution ${\cal{V}}^{\prime}$ and we get 
$$ (\nabla_XJ)(\nabla_VJ)W=0$$
for $V,W$ in ${\cal{V}}^{\prime}$ and $X$ in $H_2$. If we set ${\cal{V}}^i=H_i \bullet H_i, i=1,2$ it follows that ${\cal{V}}^1 \bullet H_2=0$. Applying the 
same reasoning to the decomposition $TM={\cal{V}}^{\prime \prime} \oplus H_1$ where ${\cal{V}}^{\prime \prime}={\cal{V}} \oplus H_1$ we get ${\cal{V}}^2 
\bullet H_1=0$. It is now routine to verify that ${\cal{V}}^1$ and ${\cal{V}}^2$ are orthogonal thus ${\cal{V}}={\cal{V}}^1 \oplus {\cal{V}}^2$. Since ${\cal{V}}^i, i=1,2$ 
are $\nabla$-parallel (one uses the formula 3.3) by the irreducibility of ${\cal{V}}$ we must have, say ${\cal{V}}^1=0$. But this would imply 
$H_1 \bullet {\cal{V}}=0$ and further $H_1 \bullet TM=0$ a contradiction, since $(M,g,J)$ is strict. q.e.d.
\\ \par
To summarize the results of this section, nearly K\"ahler manifolds with special algebraic torsion have the property that up to Riemannian products of the total 
space, the fiber (a compact, simply connected, Hermitian symmetric space) and the base space of the canonical fibration are irreducible Riemannian manifolds. 
\section{Metric properties}
In this section we will consider $(M^{2n},g,J)$, a nearly K\"ahler manifold with special algebraic torsion and 
decomposition $TM={\cal{V}} \oplus H$. We also suppose that the fiber of the canonical fibration is irreducible and recall that the base space has to be 
irreducible, too. Let us define a tensor $F : H \to H$ by 
$$F=-\sum \limits_{e_i \in {\cal{V}}}^{} (\nabla_{e_i}J)^2$$ 
where $\{ e_i \}$ is a local orthonormal frame in ${\cal{V}}$.
Our main objective will be to determine the maximal number of eigenvalues of the tensors $F$ and $C$. This will be used 
for the final classification result in the next section. It will also separate the nearly K\"ahler manifolds with 
special algebraic torsion into two classes.  \par
The fact that $C$ is $\overline{\nabla}$-parallel and symmetric, together with the irreducibilty of the fiber imply that 
there exists a real constant $\lambda$ such that $C=\lambda 1_{{\cal{V}}}$ on ${\cal{V}}$. This leads to the following: 
\begin{pro}
The tensor $C$ has at most three eigenvalues.
\end{pro}
{\bf{Proof}} : \\
We know that $-C(V \bullet X)=V \bullet CX+
CV \bullet X$ for all $V$ in ${\cal{V}}$ and $X$ in $H$ (see section 2). As ${\cal{V}}$ is an eigenspace of $C$ with eigenvalue $\lambda$ this implies that 
$S(V \bullet X)+V \bullet SX=0$ whenever $V,X$ are in ${\cal{V}}$ and $H$ respectively, where $S=C+\frac{\lambda}{2}$
on $H$. Let us denote by 
${\cal{L}}_V^H$ the projection of the Lie derivative ${\cal{L}}_V$ on $H$. Taking into account that $S^2$ is $\overline{\nabla}$-parallel  it follows using the above algebraic property of 
$S$ that ${\cal{L}}_V^H S^2=0$ for all $V$ in ${\cal{V}}$. It is now an elementary exercise to see that $S^2$ projects on a symmetric, $\nabla^h$-parallel endomorphism of $N$ which has to be 
a multiple of identity by the irreducibility of $N$, hence proving our assertion.  q.e.d.
 \\ \par 
Using the fact that $H \bullet H={\cal{V}}$ it is straightforward to show that $r_{\vert H}=2F$. Hence when 
$C$ has a single eigenvalue, that is $C=0$ it follows also that $F$ has a single eigenvalue. We concentrate now 
on the cases when $C$ has two or three eigenvalues. We will need several preliminary results.
\begin{lema}
We have : \\
(i) \begin{nr} \hfill 
$Ric((\nabla_XJ)Y,V)=\frac{1}{4}r(V, (\nabla_XJ)Y) +<Y,F(V \bullet X)+V \bullet FX>\hfill $
\end{nr}
whenever $X,Y$ are in $H$ and $V$ is in ${\cal{V}}$. \\
(ii) The base manifold $(N,h)$ is Einstein with Einstein constant $\mu >0$ and furthermore 
\begin{nr} \hfill 
$ Ric+\frac{r}{4}=\mu \cdot 1_H. \hfill $
\end{nr}
on $H$ .
\end{lema}
{\bf{Proof}} : \\
(i) is an immediate consequence of (4.2). To prove (ii), let $X,Y$ be vector fields on $N$. Then 
$Ric^h(X,Y) \circ \pi=\sum \limits_{e_ i \in H}{\cal{R}}(X^{\star}, e_i, Y^{\star},e_i)=\frac{1}{2}<r(X^{\star}), Y^{\star}>+
Ric_{\overline{\nabla}}(X^{\star},
Y^{\star})$ where $\{e_i\}$ is a local orthonormal frame in $H$. But $Ric_{\overline{\nabla}}=\frac{1}{4}(3Ric+
Ric^{\star})$ (see \cite{Gray1}) thus 
$$ Ric^h(X, Y) \circ \pi=<(Ric+\frac{r}{4})X^{\star}, Y^{\star}>.$$
Recall that $r$ is a $\overline{\nabla}$-parallel tensor, as well as $Ric$ (see \cite{Nagy1}). Hence, it is straightforward 
to see that $Ric^{h}$ is $\nabla^h$-parallel and since $(N,h)$ is irreducible, we must have $Ric^h=\mu \cdot h$ for 
some constant $\mu$. That $\mu >0$ follows by the fact that the Ricci curvature of $g$ is strictly positive (see \cite{Nagy1}). q.e.d.
\\ \par 
We mention now without proof the following elementary result.
\begin{lema}
Let $Q$ be a $\overline{\nabla}$-parallel subbundle of $H$ such that ${\cal{V}} \bullet Q=0$. Then there 
exists a $\nabla^h$-parallel subbundle $E$ of $TN$ such that $Q=\pi^{\star} E$.
\end{lema}
We are now able to treat the cases when the tensor $C$ admits exactly two or three eigenvalues. 
\begin{pro}
(i) If the tensor $C$ has exactly two eigenvalues then $F=k \cdot 1_H$ with $k>0$ constant. Furthermore, 
the eigenvalues and eigenbundles of each of the tensors $r,C,Ric$ are given in the following table. 
\begin{nr} \hfill 
$ \begin{tabular}{ | r |c | c | c |c |} 
\hline
Eigenvalue  &  r &  Ric &  C & Eigenbundle \\
\hline 
$\lambda_1 $ & $\frac{n-d}{d}k$&  $\frac{n+7d}{4d} k $&  $\frac{4(n-3d)}{d}k$ & ${\cal{V}}$ \\
\hline
$\lambda_2 $ & $2k$ & $\frac{n+2d}{2d}k$& $-\frac{2(n-3d)}{d}k$ & $H$\\
\hline
\end{tabular} \hfill $
\end{nr}
Here $d$ is the half dimension of the distribution ${\cal{V}}$. \\
(ii) Suppose that the tensor $C$ has exactly three eigenvalues. Then the tangent bundle of $M$ admits
a $\overline{\nabla}$-parallel decomposition 
$$TM={\cal{V}} \oplus H_1 \oplus H_2$$ 
with $H_1 \bullet H_1=H_2 \bullet H_2=0$ and $H_1 \bullet H_2={\cal{V}}$. Moreover, $(M,g)$ is a 
homogenous space and the base space of the canonical fibration is a simply connected, compact and 
irreducible symmetric space.
\end{pro}
{\bf{Proof}} : \\
(i)Using the (2.2) it easy to see that the eigenvalues of $C$ are of the form 
$\lambda$ and $-\frac{\lambda}{2}$. Since by definition $C=Ric-5Ric^{\star}$ we get by (5.2) that 
$r=\frac{-\lambda+8\mu}{12}$ on $H$ hence $F$ has a single eigenvalue.\\
(ii) Let $\lambda_1$ and $\lambda_2$ be the eigenvalues of $C$ on $H$, and denote the corresponding 
eigenbundles by $H_1$ and $H_2$. Using lemma 5.2 we have that ${\cal{V}} \bullet H_1 \neq 0$, hence it follows 
by (2.2) that $-(\lambda+\lambda_1)$ is an eigenvalue of $C$. It is different from $\lambda$ since ${\cal{V}} \bullet 
H_1$ is orthogonal to ${\cal{V}}$. Let us suppose that $-(\lambda+\lambda_1)=\lambda_1$. 
It follows that ${\cal{V}} \bullet H_1 \subseteq H_1$ and since $M$ has special algebraic torsion we get that $H_2 \bullet 
{\cal{V}}$ is orthogonal ot ${\cal{V}} \oplus H_1$, that is $H_2 \bullet 
{\cal{V}}$ is contained in $H_2$. Again by lemma 5.2 we have $H_2 \bullet {\cal{V}} \neq 0$ and using (2.2) we 
find that $\lambda_2=-\frac{\lambda}{2}=\lambda_1$ an absurdity. \par
Hence $-(\lambda+\lambda_1)=\lambda_2$ a fact which obviously implies by (2.2) that ${\cal{V}} \bullet 
H_1 \subseteq H_2, {\cal{V}} \bullet H_2 \subseteq H_1$ and $H_1\bullet H_2 \subseteq {\cal{V}}$. Let us show now 
that $H_1 \bullet H_1=0$. If $H_1 \bullet H_1 \neq 0$ it produces the eigenvalue $-2\lambda_2$ for $C$. Or 
$H_1 \bullet H_1$   is orthogonal to $H_2$ and $-2\lambda_1 \neq \lambda$ (since $\lambda_1 \neq \lambda_2$). We 
get $-2\lambda_1=\lambda_1$ hence $\lambda_1=0$ and $H_1 \bullet H_1 \subseteq H_1$. We set 
${\cal{V}}^{\prime}=H_1$ and $H^{\prime}={\cal{V}} \oplus H_2$. Then using corollary 3.2, (ii) we obtain that
$(\nabla_XJ)(\nabla_YJ)Z$ belongs to $H^{\prime}$ whenever $X,Y$ are in ${\cal{V}}$ and $Z$ is in $H_1 \bullet 
H_1$. As ${\cal{V}} \bullet ( {\cal{V}} \bullet H_1) \subseteq {\cal{V}} \bullet H_2 \subseteq H_1$ we obtain 
easily that ${\cal{V}} \bullet (H_1 \bullet H_1)=0$ and we conclude by lemma 5.2. In the same way it can be proven 
that $H_2 \bullet H_2=0$. Hence, we have three different splittings of special algebraic type 
$TM={\cal{V}} \oplus H, TM=H_1 \oplus ({\cal{V}} \oplus H_2), TM=H_2 \oplus 
({\cal{V}} \oplus H_1)$ each of which being $\overline{\nabla}$-parallel. It follows that the restriction of the tensor 
$\overline{\nabla}_U \overline{R}$ to either ${\cal{V}}, H_1$ or $H_2$ vanishes. Or this implies easily that 
$\overline{\nabla} \ \overline{R}=0$, so $\overline{\nabla}$ is an Ambrose-Singer connection hence $(M,g,J)$ 
is a homogenous space \cite{Tri}. That the base space of the canonical fibration is symmetric follows by the 
usual comparison between curvature tensors of the total and the base space in \cite{Besse1}. q.e.d.
$\\$ \par
Hence we obtained a new class of nearly K\"ahler manifolds motivating the following definition.
\begin{defi}
A nearly K\"ahler manifold is homogenous of type $III$ if it has special algebraic 
torsion and the tensor $C$ admits exactly three eigenvalues.
\end{defi}
\begin{rema}
In the case when $C$ has exactly three eigenvalues, let $d_1$ and $d_2$ be the 
half dimensions of $H_1$ and $H_2$ respectively. 
Then the eigenvalues of the tensors 
$r, Ric, C$ with their corresponding eigenbundles are listed in the following table.
\begin{nr} \hfill
$\begin{tabular}{ | r |c | c | c |c | }
\hline 
Eigenvalue  &  r &  Ric &  C & Eigenbundle \\
\hline
$\lambda_1 $ & $\frac{2d_1}{d}k$&  $(\frac{d_1}{2d}+\frac{d_1}{d_2}+1)k
$&  $4(\frac{2d_1}{d}-\frac{d_1}{d_2}-1)$ & ${\cal{V}}$ \\
\hline
$\lambda_2 $ & $k$ & $(\frac{d_1}{d}+\frac{d_1}{d_2}+\frac{1}{2})k
$& $-4(\frac{d_1}{d}+\frac{d_1}{d_2}-2)$& $H_1$\\
\hline
$\lambda_3 $ & $\frac{d_1}{d_2}k$ & $(\frac{d_1}{d}+\frac{d_1}{2d_2}+1)k
$& $-4(\frac{d_1}{d}-\frac{2d_1}{d_2}+1)$ & $H_2$\\
\hline
\end{tabular} \hfill $
\end{nr}
Here, $k$ is a positive constant and, moreover, the eigenvalues 
of the tensor $F$ are $k$ and $\frac{d_1}{d_2}k$ with eigenbundles $H_1$ and $H_2$ respectively. The proof will be omitted, as being a simple calculation, 
based on (2.3) and of the relations between eigenvalues of the tensors $C,F$ and $r$ developed in the proof of the proposition 5.2, (ii).
\end{rema}
To summarize the results obtained so far in this section, one can reduce the study of special algebraic torsion to the case when 
the tensor $F$ has a single eigenvalue. 
\section{The holonomy of the base manifold}
We consider in this section a nearly K\"ahler manifold with special algebraic torsion and let $TM={\cal{V}} \oplus H$ 
be the corresponding decomposition. Furthermore, we suppose that the fiber of the canonical fibration is irreducible 
and that the tensor $F$ has a single eigenvalue, to be denoted by $k$. \par
Let $m$ be a point of $M$. We define a vector subspace of 
$\mathfrak{so}(H_m)$ by 
$$ \mathfrak{p}_m=\{ \nabla_vJ : v \ \mbox{in} \ {\cal{V}}_m \}.$$
The dimension of $\mathfrak{p}_m$ equals that of ${\cal{V}}_m$ since $(M,g,J)$ is strict. Let $\mathfrak{k}_m$ be the vector 
subspace of $\mathfrak{so}(H_m)$ generated by 
$$ \{ [ \nabla_vJ, \nabla_wJ] : v,w \in {\cal{V}}_m \}. $$
\begin{pro}
We have $ \mathfrak{k}_m \cap \mathfrak{p}_m=0 $. Then $ \mathfrak{h}_m= \mathfrak{k}_m \oplus \mathfrak{p}_m $ is a Lie 
subalgebra of $\mathfrak{so}(H_m).$
\end{pro}
{\bf{Proof}} : \\
Let us first prove that $[\mathfrak{p}_m, \mathfrak{k}_m] \subseteq \mathfrak{p}_m$. Indeed, if $V_1, V_2, V_3$ belong to 
${\cal{V}}$ we obtain using (4.2) that 
$$[\nabla_{V_1}J, [\nabla_{V_2}J,\nabla_{V_3}J]=-\nabla_{R(V_2,V_3)V_1}J $$
and the assertion is proven. Consider now $z$ in $\mathfrak{p}_m \cap \mathfrak{k}_m$. Then $z=\nabla_vJ$ with 
$v$ in ${\cal{V}}$ and $[\nabla_vJ, x]$ belongs to $\mathfrak{p}_m$ for all $x$ in $\mathfrak{p}$. In particular, there 
exists $w$ in ${\cal{V}}$ such that $[\nabla_vJ, \nabla_{Jv}J]=\nabla_wJ$, which implies $2(\nabla_vJ)^2=
\nabla_{Jw}J$. The left side of this equality is symmetric whilst the right is skew-symmetric. It follows that $z=0$ 
and we proved that $ \mathfrak{k}_m \cap \mathfrak{p}_m=0 $. Now, $[\mathfrak{p}_m,\mathfrak{p}_m] \subseteq 
\mathfrak{k}_m$ by definition and using the Jacobi identity and (4.2) we obtain 
$$ [[\nabla_{V_1}J, \nabla_{V_2}J],[\nabla_{V_3}J,\nabla_{V_4}J]]=
[\nabla_{V_2}J, \nabla_{R(V_3,V_4)V_1}J]+[\nabla_{R(V_3,V_4)V_2}J, \nabla_{V_1}J]$$
for all $V_i$ in ${\cal{V}}, 1 \le i \le 4$. It follows that $\mathfrak{k}_m$ and $\mathfrak{h}_m$ are Lie algebras. q.e.d.
\begin{rema}
(i) We saw along the proof of the previous proposition that $[\mathfrak{p}_m, \mathfrak{p}_m] \subseteq 
\mathfrak{k}_m$ and $[\mathfrak{p}_m, \mathfrak{k}_m] \subseteq \mathfrak{p}_m$. \\
(ii) Using E. Cartan's description of symmetric spaces one can also relate the Lie algebras $\mathfrak{h}_m, \mathfrak{k}_m$ to the homogenenous 
structure of the symmetric space $F_m$, the fiber of $\pi$ through $\pi(m)$.
\end{rema}
We will now relate the Lie algebra $\mathfrak{h}_m$ to the holonomy group of the base manifold. We define the Lie 
algebra $\mathfrak{hol}_m$ as the Lie subalgebra of $\mathfrak{so}(H_m)$ generated by skewsymmetric endomorphisms of the 
form 
$$\tau^{-1}_c \circ {\cal{R}}(\tau_cv, \tau_cw) \circ \tau  $$
where $v,w$ belong to $H_m$, $c$ is a horizontal curve in $M$ starting at $m$ and $\tau_c$ is parallel transport along 
$c$ with respect to the connection $\overline{\nabla}$. Using standard considerations in holonomy theory we find that 
$\mathfrak{hol}_m$ is in fact isomorphic with the Lie algebra of the holonomy group of $(N,h)$ through $\pi m$. 
\begin{lema}
The Lie algebra $\mathfrak{h}_m$ is an ideal of $\mathfrak{hol}_m$.
\end{lema}
{\bf{Proof}} : \\
Let us first note that (3.3) gives us, after some computations 
\begin{nr} \hfill 
$ [\overline{R}(X,Y), \nabla_VJ]=\nabla_{\overline{R}(X,Y)V}J \hfill $
\end{nr}
for all $X,Y$ in $H$ and $V$ in ${\cal{V}}$. It follows that $[{\cal{R}}(v,w), \mathfrak{p}_m] \subseteq 
\mathfrak{h}_m$ and since $[\mathfrak{p}_m, \mathfrak{p}_m]=\mathfrak{k}_m$ by the Jacobi 
identity we also get $[{\cal{R}}(v,w), \mathfrak{k}_m] \subseteq 
\mathfrak{h}_m$ for all $v,w$ in $H_m$. Let $c$ be a horizontal path in $M$ starting at $m$. Using (3.3) we obtain easily that 
$ \tau_c \circ \nabla_vJ=(\nabla_{\tau_c v}J) \circ \tau_c$ whenever 
$v$ belongs to $H_m$. The conclusion now follows. q.e.d.
\\ \par
Let us consider now the bundle $\mathfrak{h}$ whose fiber at a point $m$ of $M$ equals $\mathfrak{h}_m$. 
Note that $\mathfrak{h}$ can be identified, in a standard way, with a subbundle of $\Lambda^2(H)$. We do this 
identification tacitly in the rest of this section. In the 
same way we obtain bundles $\mathfrak{k}, \mathfrak{p}$. We need now the following metric fact.
\begin{pro}
We have $\rho^{{\cal{R}}}=\frac{k(n-d)}{d}$ on $\mathfrak{h}$ where we denoted by 
$\rho^{{\cal{R}}}$ the curvature operator of ${\cal{R}}$. 
\end{pro}
{\bf{Proof}} : \\
That the claimed formula holds on $\mathfrak{p}$ was proved in lemma 4.2 (see also table 5.3). Let us show that 
it holds on $\mathfrak{k}$. We must compute $\sum \limits_{e_i \in H}^{} {\cal{R}}(e_i,qe_i,X,Y)$ where $e_i$ is a local orthonormal basis in $H$ and $q=[\nabla_VJ,
\nabla_WJ]$ with $V,W$ in ${\cal{V}}$. The definition of $C$ implies easily that our quantity equals $\sum \limits_{e_i \in H}^{} {\overline{R}}(e_i,qe_i,X,Y)+
\sum \limits_{e_i \in H}^{} \nabla_{(\nabla_{e_i}J)qe_i}J$. Or $\sum \limits_{e_i \in H}^{} (\nabla_{e_i}J)qe_i=0$ (one uses the base ${Je_i}$). Now, using the first Bianchi 
identity for the Hermitian connection $\overline{\nabla}$ (see \cite{KN}) we get  : 
$$ \begin{array}{lr}
\sum \limits_{e_i \in H}^{}\overline{R}(e_i, qe_i, X, Y)+\sum \limits_{e_i \in H}^{}\overline{R}(qe_i, X,e_i,Y)+
\sum \limits_{e_i \in H}^{}\overline{R}(X,e_i,qe_i,Y)=\vspace{2mm}\\
=Tr_H [\nabla_YJ, \nabla_XJ]q. 
\end{array} $$
But $\sum \limits_{e_i \in H}^{}\overline{R}(qe_i, X,e_i,Y)=\sum \limits_{e_i \in H}^{}<[\overline{R}(X,e_i),q]e_i,Y>+
Ric_{\overline{\nabla}}(X,qY)$ where $Ric_{\overline{\nabla}}$ denotes the Ricci tensor of the Hermitian connection 
$\overline{\nabla}$. The use of (6.1) implies that 
$$<[\overline{R}(X,e_i),q]e_i,Y>=
<[\nabla_VJ,\nabla_{\overline{R}(X,e_i)W}J]e_i,Y>-<[\nabla_WJ,\nabla_{\overline{R}(X,e_i)V}J]e_i,Y>.$$
Now using lemma 3.1 we obtain after computing at some length that 
$$\begin{array}{lr}
 \sum \limits_{e_i \in H}^{}<[\nabla_VJ,\nabla_{\overline{R}(X,e_i)W}J]e_i,Y>=-<(\nabla_WJ)X, G(\nabla_VJ)Y)>-\\
-Tr_H(\nabla_XJ)(\nabla_YJ)(\nabla_VJ)(\nabla_WJ).
\end{array}$$
where $G: H \to H$ is defined by $G=-\sum \limits_{e_i \in H} (\nabla_{e_i}J)^2.$ It follows that 
$$\begin{array}{lr}
\sum \limits_{e_i \in H}^{}<[\overline{R}(X,e_i),q]e_i,Y>=-Tr_H(\nabla_XJ)(\nabla_YJ)q+
<(\nabla_VJ)X,G(\nabla_WJ)Y>-\\
<(\nabla_WJ)X, G(\nabla_VJ)Y>.
\end{array}$$
Since $G=k1_H$ we finally obtain that 
$$\sum \limits_{e_i \in H}^{}\overline{R}(e_i, qe_i, X,Y)+Ric_{\overline{\nabla}}(X,qY)-Ric_{\overline{\nabla}}(Y,qX)
+2k<qX,Y>=0. $$
We conclude by the fact (see \cite{Gray1}) that $Ric_{\overline{\nabla}}=\frac{1}{4}(3Ric+Ric^{\star})$ on $M$ and 
by using the proposition 5.2, (i). q.e.d.
\\ \par
We are now ready to prove : 
\begin{teo}
Let $(M,g,J)$ be a nearly K\"ahler manifold with special algebraic torsion and let us suppose that the fiber of the 
canonical fibration is irreducible. If the tensor $F$ has a single eigenvalue and 
the base space of the canonical fibration is not symmetric then $M$ is the twistor 
space over a quaternionic-K\"ahler manifold of positive scalar curvature.
\end{teo}
{\bf{Proof}} : \\
The base space $(N^{2m},h)$ of the canonical fibration is an irreducible, Einstein manifold of strictly 
positive scalar curvature. Then using the Berger holonomy theorem \cite{pitman} and the well known fact 
that metrics with holonomy $SU(m), Sp(q), G_2$ or $Spin(7)$ are Ricci flat we find that that they are only 
three posibilities for the holonomy group of $(N,h)$, namely $SO(2m), U(m)$ and $Sp(q) \cdot Sp(1)$ when 
$m=2q$. We will now treat each case separately. \par
If $Hol(N,h)=SO(2m)$, then if $m \neq 3$ the Lie algebra $\mathfrak{so}(2m)$ is simple and by lemma 6.1 we get easily that $\mathfrak{h}=
\Lambda^2(H)$ so $(N,h)$ has constant sectional curvature by proposition 6.2. Hence $(N,h)$ is a symmetric space 
(in fact a round sphere), an absurdity. If $m=2$ then $\mathfrak{so}(4)$ has two simple factors, each of dimension 
$3$. Thus if $\mathfrak{h} \neq \Lambda^2(H)$ it has dimension $3$ and it follows that ${\cal{V}}$ has dimension 
$2$. We found that $(M,g,J)$ is $6$-dimensional, Hermitian reducible nearly K\"ahler manifold, hence $(N,h)$ 
is symmetric by \cite{Moro1}. \par
We consider now the case when $(N,h)$ has holonomy 
equal to $U(m)$. Then $(N,h)$ is a K\"ahler manifold with associated complex structure $I$. The holonomy bundle of 
$(N,h)$ is then $\Lambda^{1,1}(N)$. Denote by $\Lambda_0^{1,1}(N)$ the space of traceless two forms of type $(1,1)$.
Using lemma 4.1 and the fact that $\mathfrak{h}$ has rank at least $2$ we have either that $\mathfrak{h}=\pi^{\star} \Lambda^{1,1}(N)$ or  
$\mathfrak{h}=\pi^{\star} \Lambda^{1,1}_0(N)$. As the curvature tensor of a K\"ahler Einstein manifold is determined 
by its restriction to $ \Lambda^{1,1}_0(N)$  it is easy to conclude using proposition 6.2 that $(N,h)$ is symmetric, 
a contradiction. \par
Lastly, we suppose that $m=2q$ and $Hol(N,h)=Sp(q) \cdot Sp(1)$. Thus, 
$(N,h)$ is a quaternionic K\"ahler manifold and let us denote by $Q$ the $3$-dimensional subbundle of 
$\Lambda^2(N)$ giving the quaternionic structure. Then we have a decomposition $\Omega^2(N)=E_1 \oplus E_2 \oplus Q$ 
where fibers of $E_1$ and $Q$ are isomorphic in each point with $\mathfrak{sp}(q)$ and $\mathfrak{sp}(1)$ respectively. Furthermore, 
the holonomy bundle of the metric $h$ equals $E_1 \oplus Q$. Now the curvature transformation of 
of the metric $h$ is determined by its restriction to $E_1$ as shown in \cite{Sal3}. Hence if it is constant on $E_1$, then $(N,h)$ is symmetric, in 
fact a quaternionic projective space. It follows that the only possibility is that $\mathfrak{h}=\pi^{\star}Q$ and this implies that ${\cal{V}}$ is of 
rank  two. Thus, the canonical Hermitian connexion of $(M,g,J)$ has holonomy contained in $U(1) \times U(2q)$ and we conclude by 
using a result from \cite{Nagy1}. q.e.d.
$\\$ \par 
It remains to investigate the case when the manifold $(N,h)$ is a symetric space. We set the following : 
\begin{defi}
A simply connected, nearly K\"ahler manifold with special algebraic torsion and such that the canonical fibration has irreducible fiber is 
called homogenous of type 
IV if the tensor $F$ has exactly two eigenvalues and the base space of the canonical fibration is a symmetric space.
\end{defi}
Let us remark that the terminology "homogenous" is not yet justified. However, using the relation between the tensors 
$\overline{R}$ and $R^h$ given by O'Neill's relations, one easily gets that when $N$ is symmetric, 
$\overline{\nabla}$ is an Ambrose-Singer connection, hence $M$ is homogenous. \par
Now the proof of the theorem 1.1 is a consequence of the material previously presented. 
\section{Nearly K\"ahler metrics from Riemannian foliations}
The main purpose of this section is to give a proof of the theorem 1.2, stated in the introduction. Let 
us consider a K\"ahler manifold $(M,g,J)$ together with a Riemannian foliation ${\cal{F}}$. Let ${\cal{V}}$ be the 
integrable distribution induced by ${\cal{F}}$. Following \cite{Tondeur} we recall that the metric $g$ is bundle-like with respect to our foliation, that is 
${\cal{L}}_Vg$ vanishes on $H$, the orthogonal complement of ${\cal{V}}$, for every $V$ in ${\cal{V}}$. As in the statement of the theorem 1.3 we make 
the assumption that the leaves of ${\cal{F}}$ are complex submanifolds of $M$, which gives that  $J {\cal{V}} \subseteq {\cal{V}}$. \par
We consider now the Riemannian metric on $M$ defined by 
$$ \hat{g}(X,Y)=\frac{1}{2}g(X,Y) \ \mbox{if} \ X,Y \in {\cal{V}},  
\hat{g}(X,Y)=g(X,Y) \ \mbox{for} \ X,Y \mbox{in} \ H.$$
The metric $\hat{g}$ admits a compatible almost complex structure 
$\hat{J}$ given by $\hat{J}_{\vert {\cal{V}}}=-J$ and $\hat{J}_{\vert H}=J$. This almost complex structure 
was introduced in \cite{Eells} for the case of twistor spaces over $4$-manifolds. 
\begin{pro}
The manifold $(M,\hat{g},\hat{J})$ is nearly K\"ahler. The distributions 
${\cal{V}}$ and $H$ are parallel with respect to the canonical Hermitian 
connection of $(M,\hat{g},\hat{J})$ which thus has reduced holonomy.
\end{pro}
{\bf{Proof}} : 
Let $A : TM \times TM \to TM$ be the O'Neill tensor of the Riemannian foliation induced by ${\cal{V}}$.  As 
$g$ is K\"ahler we must have $A_XJ=JA_X$ for all $X$ in $TM$. Using the relations 
between the Levi-Civita connections of $\hat{g}$ and $g$ given in \cite{Tondeur} we obtain 
after a standard computation : 
$$ \begin{array}{lr}
(\hat{\nabla}_X \hat{J})V=-(\hat{\nabla}_V \hat{J})X=-A_X(JV) \\
(\hat{\nabla}_V \hat{J})W=0, \ (\hat{\nabla}_X \hat{J})Y=2A_X(JY)
\end{array} $$
for every $X,Y$ in ${\cal{V}}$ and $V,W$ in $H$. It is now straightforward to conclude. q.e.d.
\begin{coro}
The twistor space of a positive quaternionic-K\"ahler manifold of 
dimension $4k$ admits a canonical  NK structure with reducible holonomy, 
contained in $U(1) \times U(2k)$.
\end{coro}
{\bf{Proof}} : \\
We have only to recall \cite{Sal3} that such a twistor space is the total space of a Riemannian submersion with 
totally geodesic fibers of dimension $2$ and that it admits a compatible K\"ahler structure. q.e.d.
\begin{rema}
If one drops the condition that the distibution ${\cal{V}}$ be totally geodesic, the construction above produces only a 
quasi-K\"ahler manifold. 
\end{rema}
We are now going to prove theorem 1.2. We remark that the nearly K\"ahler manifold $(M,\hat{g}, \hat{J})$ has the property that ${\cal{V}} \bullet {\cal{V}}=0$ and 
$H \bullet H \subseteq {\cal{V}}$. Using the decomposition result  of the proposition 2.1 we obtain that (under the simple connectivity assumption) $M$ is the 
Riemannian product of a K\"ahler manifold and a nearly K\"ahler manifold with special algebraic torsion. Hence by the structure results of sections 4,5 and 6 the de deRham decomposition 
of the second factor consists only of spaces of type III and IV and twistor spaces over positive quaternionic K\"ahler manifolds equipped with their canonical nearly K\"ahler metric. The proof 
of the theorem comes now by reversing the construction done in the beginning of the section, and by using proposition 4.3. However, we postpone until the next section the 
proof of the fact that a homogeneous space $(M,g,J)$ of type III or IV is a homogeneous K\"ahler manifold  when endowed with the metric $\overline{g}$ and the complex structure 
$\overline{J}$.
\section{Homogeneous nearly K\"ahler manifolds}
In this section we will investigate some properties of homogeneous nearly K\"ahler manifolds that are introduced by the following definition. 
\begin{defi}
Let $(M^{2n},g,J)$ a strict and complete nearly K\"ahler manifold. It is called homogeneous if it admits a transitive action of a connected, closed 
subgroup of the group of holomorphic isometries of $(g,J)$.
\end{defi}
Let us now prove that what we called homogenous spaces of type I, II, III, IV are in fact homogeneous nearly K\"ahler manifolds and that for the spaces of type 
III and IV the natural K\"ahler metric is homogeneous.
\begin{pro}
Let $(M^{2n},g,J)$ be a strict and complete nearly K\"ahler manifold. \\
(i) If $\overline{\nabla} \ \overline{R}=0$ then $(M,g,J)$ is a homogeneous 
nearly K\"ahler space. In particular $M$ is naturally reductive and if moreover, $M$ is simply connected, then $(M,g,J)$ is a Riemannian $3$-symmetric space. \\
(ii) If $(M,g,J)$ is a space of type III or IV then $(M, \overline{g}, \overline{J})$ is a homogeneous K\"ahler manifold.
\end{pro} 
{\bf{Proof}} : \\
(i) In the terminology of \cite{Tri}, the $\overline{\nabla}$-parallel tensor $T_XY=-\frac{1}{2}(\nabla_XJ)JY$ is a homogeneous structure. The theorem of Ambrose and Singer (see \cite{Tri})
gives now an explicit construction of $M$ as a homogeneous space and from this it is easy to conclude using a result of Gray (see \cite{Gray4}, proposition 5.5, page 358). \\
(ii) Using the relation between the Levi-Civita of $\overline{g}$ and $\overline{\nabla}$ given in the proof of proposition 4.3 one can prove 
directly that the tensor $T_XY$ below is a homogeneous structure for $(M,\overline{g}, \overline{J})$. We omit the details. q.e.d.
\\ \par
Let us give now a corollary of theorem 1.1 related to 
homogeneous nearly K\"ahler spaces.
\begin{pro}
Let $(M^{2n},g,J)$ be a homogeneous nearly K\"ahler manifold which is also strict and simply connected. Then $(M,g,J)$ is isometric and 
biholomorphic to a Riemannian product whose factors 
belongs to the following list : \\
-homogeneous $6$-dimensional \\
-spaces of type I,II,III or IV \\
-twistor spaces over positive quaternionic K\"ahler manifolds which are homogeneous as nearly K\"ahler manifolds.
\end{pro}
{\bf{Proof}} : \\
It is not hard to see that each of the manifolds occuring in the decomposition of theorem 1.1 is irreducible in the usual Riemannian sense. Hence, that decomposition 
is the deRham decomposition of the nearly K\"ahler manifold $(M,g,J)$ and by a result of Hano \cite{KN} it follows that the connected isometry group of $(M,g)$ is the direct 
product of the isometry groups of the factors. But this implies that each factor is a homogeneous space and using that $J$ preserves our decomposition it is straightforward 
to conclude q.e.d.
\\ \par
The rest of this section will be devoted to the proof of the following classification result. 
\begin{teo}
Let $(M^{2n},g,J)$ be a homogeneous nearly K\"ahler manifold which is also strict and simply connected. Then $(M,g,J)$ is isometric and bihomorphic to a Riemannian product whose factors 
belong to the following list : \\
-homogeneous $6$-dimensional nearly K\"ahler manifolds \\
-spaces of type I,II,III or IV \\
\end{teo}
\par Before getting into the proof let us make some comments on theorem 8.1. As any Riemannian 
$3$-symmetric space is nearly K\"ahler as soon as it is naturally reductive the previous result yields a 
geometric version of the classification of naturally reductive 
$3$-symmetric spaces. Note that the general classification theory of $3$-symmetric spaces develloped 
in \cite{Gray4} was obtained using quite involved Lie algebra arguments. It also follows from theorem 8.1 
that Wolf\&Gray's conjecture (see \cite{Wolf2}, page 158) asserting that any strict homogenenous 
nearly K\"ahler manifold is $3$-symmetric holds true if and only if it holds in $6$-dimensions. Concerning the types 
of homogenenous spaces  appearing in the theorem let us recall that those of type I are compact, simply connected 
isotropy irreducible homogeneous spaces, of naturally reductive type, carrying a nearly K\"ahler structure. Here 
our geometric results do not allow classification. For spaces of type II, further work is required 
to localize completely this class in Gray's list (see section 3, remark 3.1). Finally, spaces of types III and IV can be classified using 
on the one hand the classification of compact and simply connected, irreducible Hermitian symmetric spaces and other 
the other hand the description of the fiber of the canonical fibration in the terms of the Lie algebra $\mathfrak{h}$ 
together with lemma 6.1 (see \cite{Nagy3} for details). Note that the spaces of type III and IV were already studied from various points of view in \cite{Sal86} and 
\cite{Burstall}. 
\par
Let us now prove the theorem 8.1. It suffices to show that if a twistor space with its canonical nearly K\"ahler metric is homogeneous then the base quaternionic manifold 
is symmetric. For, let us consider such a twistor space; it is a strict nearly K\"ahler manifold with special algebraic torsion $(M,g,J)$  with decomposition  
$TM={\cal{V}} \oplus H$ where ${\cal{V}}$ is of rank two and such the tensor $F$ has a single eigenvalue of the form $\frac{k}{4}, k >0$. \par 
Let us set now the following notation. If $(Z,h)$ is a Riemannian manifold we denote by $\iota(Z,h)$ the Lie algebra of Killing vector 
fields. Furthermore, if $I$ is an almost complex structure giving $(Z,h)$ the structure of a almost complex manifold we define 
$\iota(Z,h,J)$ to be the Lie algebra of holomorphic Killing vector vector fields. \par 
Our aim is to obtain a comparison result between $\iota(M,\overline{g},\overline{J})$ and $\iota(M,g,J)$. Let 
$X=V+K$ be an element of $\iota(M,g,J)$ split into its horizontal and vertical components. We use the parallelism of the decomposition 
$TM={\cal{V}} \oplus H$ with respect to the connection $\overline{\nabla}$ in order to see that the Killing equation for $X$ 
(the skew-symmetry of $\nabla X$) is equivalent with the following system of equations for its components : 
\begin{nr} \hfill 
$ \begin{array}{lr}
<\overline{\nabla}_YK, Z>+ <\overline{\nabla}_ZK,Y>=0 \hspace{0.5cm} \\
<\overline{\nabla}_WK, Y>+ <\overline{\nabla}_YV,W>=0 \hspace{0.5cm} \\
<\overline{\nabla}_{W_1}V, W_2>+ <\overline{\nabla}_{W_2}W_1>=0 \\
\end{array} \hfill $ 
\end{nr}
where $Y,Z$ are in $H$ and $W,W_1,W_2$ are in ${\cal{V}}$. Using the same arguments the condition  that $X$ be a $J$-holomorphic 
vector field projets onto the following two systems : 
\begin{nr} \hfill 
$\begin{array}{lr}
\overline{\nabla}_{JY}K-J(\overline{\nabla}_YK)=-2(\nabla_YJ)V 
\hspace{0.5cm} \\
\overline{\nabla}_{JY}V-J(\overline{\nabla}_YV)=-2(\nabla_YJ)K 
\end{array} \hfill $ 
\end{nr}
and 
\begin{nr} \hfill 
$ \begin{array}{lr}
\overline{\nabla}_{JW}K-J(\overline{\nabla}_WK)=-2(\nabla_WJ)K 
\hspace{0.5cm} \\
\overline{\nabla}_{JW}V-J(\overline{\nabla}_WV)=0 
\end{array} \hfill $ 
\end{nr}
whenever $Y$ and $W$ are in $H$ and ${\cal{V}}$ respectively. 
\begin{rema}
Neither all the equations of the systems (8.1), (8.2), (8.3) are independent, nor all will be used in 
what follows. We gave all of them for completeness reasons.
\end{rema}
We will need now the following auxiliary result. 
\begin{lema}
If $Y$ and $W$ are in $H$ and ${\cal{V}}$ respectively, then : \\ 
(i) $\overline{\nabla}_WV=\frac{f}{2}JW$ where $f=d^{\star}(JV)$, \\
(ii) $\overline{\nabla}_WK=(\nabla_WJ)JK$, \\
(iii) $\overline{\nabla}_Y V=(\nabla_YJ)JK$.
\end{lema}
{\bf{Proof}} : \\
The proof of (i), an easy exercise, will be left to the reader. In order to prove (ii) let $\{ e_i\}, i=1,2$ be a local orthonormal basis 
in ${\cal{V}}$. Using the second equation of (8.1) with $W=e_i$ and 
next deriving it in the direction of $e_i$ we obtain : 
$$ <\overline{\nabla}_{e_i}\overline{\nabla}_{e_i}K,Y>+<\overline{\nabla_{e_i}}K, \overline{\nabla}_{e_i}Y>+
<\overline{\nabla}_{e_i}\overline{\nabla}_{Y}V,e_i>+<\overline{\nabla}_{Y}V, \overline{\nabla}_{e_i}e_i>=0.$$
Using again the second equation of (8.1) for the second and third terms of the previous sum we get that $<\overline{\nabla}^2_{e_i,e_i}K,Y>+
<\overline{\nabla}_{e_i,Y}^2V,e_i>=0$. But by the Ricci identity $\overline{\nabla}_{e_i,Y}^2V=\overline{\nabla}_{Y, e_i}^2V-\overline{R}(e_i,Y)V+
\overline{\nabla}_{(\nabla_YJ)Je_i}V$ and moreover $\overline{R}(e_i,Y)V=0$ as $e_i$ is in ${\cal{V}}$ and $Y$ belongs to $H$. Now , 
$$ <\overline{\nabla}_{Y, e_i}^2V, e_i>=<\overline{\nabla}_Y\overline{\nabla}_{e_i}V-\overline{\nabla}_{\overline{\nabla}_Ye_i}V,e_i>=
\frac{1}{2}Y.f <Je_i,e_i>=0$$
by (i). Hence $<\overline{\nabla}^2_{e_i,e_i}K,Y>+<\overline{\nabla}_{(\nabla_YJ)Je_i}V,e_i>=0$ and using  the second equation of (8.1) 
(together 
with the fact that $\nabla_{e_i}J$ is skew-symmetric and $J$-anticommuting) we arrive at $<\overline{\nabla}^2_{e_i,e_i}K,Y>=
<(\nabla_{e_i}J)J \overline{\nabla}_{e_i}K, Y>$. We define now 
$$ S= \sum \limits_{i=1}^{2} (\nabla_{e_i}J)J \overline{\nabla}_{e_i}K.$$
Then $S=\sum \limits_{i=1}^{2} (\nabla_{Je_i}J)J \overline{\nabla}_{Je_i}K$ and we use that $(\nabla_{JU}J)JT=-(\nabla_{U}J)T$ for all 
$U,T$ in $TM$ and also the first equation of (8.3) to show that $S=\sum \limits_{i=1}^{2} (\nabla_{e_i}J)^2K=-\frac{k}{2}K$. Resuming, 
we have that $\sum \limits_{i=1}^{2}<\overline{\nabla}^2_{e_i,e_i}K,Y>=-\frac{k}{2}<K,Y>$. On the other hand, as 
$\overline{\nabla}_{e_i}e_i=
\nabla_{e_i}e_i$ and by the second equation of the system (8.1) we obtain that 
$$ 
\begin{array}{lr}
\Vert \overline{\nabla}_{e_i}K\Vert^2=e_i.<\overline{\nabla}_{e_i}K,K>-<\overline{\nabla}_{\nabla_{e_i}e_i}K,K>-
<\overline{\nabla}^2_{e_i,e_i}K,K>=\\
-\biggl [ e_i.<\overline{\nabla}_KV,e_i>-<\overline{\nabla}_KV, 
\nabla_{e_i}e_i>\biggr ]-<\overline{\nabla}_{e_i, e_i}^2K,K>.
\end{array} $$ 
As the vector field $\overline{\nabla}_KV$ is vertical, we obtain easily that 
$ \sum \limits_{i=1}^{2} \Vert \overline{\nabla}_{e_i}K \Vert^2=\frac{k}{2}\Vert K\Vert^2+d^{\star} \alpha $
where $\alpha$ is the $1$-form of $M$ dual to $\overline{\nabla}_KV$. By the previous considerations, this implies 
that 
$$ \sum \limits_{i=1}^{2} \Vert \overline{\nabla}_{e_i}K-(\nabla_{e_i}J)JK \Vert^2=d^{\star} \alpha $$
and the conclusion follows by integration over $M$. The identity in (iii) is a consequence of (i) and of 
the second equation of (8.1). q.e.d.
$\\$
\par
We can now prove the following.
\begin{pro}
We have $\iota(M,g,J)=\iota(M,\overline{g}, \overline{J})$. 
\end{pro}
{\bf{Proof}} : \\
Let us prove that $\iota(M,g,J) \subseteq \iota(M,\overline{g}, \overline{J})$. If $\widetilde{\nabla}$-denotes 
the Levi-Civita connection of $(M,\overline{g})$ then we recall that 
$\widehat{\nabla}_XY=\nabla_XY, \widetilde{\nabla}_XV=\overline{\nabla}_XV-(\nabla_XJ)JV, 
\widetilde{\nabla}_VX=\overline{\nabla}_VX $ and $\widetilde{\nabla}_VW=\nabla_VW$ for all $X,Y$ in $H$ and 
$V,W$ in ${\cal{V}}$ respectively. For, if $X$ in $\iota(M,g,J)$ it is easy to see, using the previous lemma that the 
equation $\widetilde{\nabla}X=0$ follows from the system 8.1 (of course we have to take into account that on 
a compact K\"ahler manifold every Killing vector field is holomorphic). The other inclusion can be proven, for example, 
using the same technique, and will be left to the reader. q.e.d.
$\\$ \par
It follows that $Isom^{0}(M,g,J)=Isom^{0}(M,\overline{g}, \overline{J})$, hence if 
$(M,g,J)$ is a homogeneous nearly K\"ahler manifold then $(M,\overline{g}, \overline{J})$ is a homogeneous 
K\"ahler manifold. Recall now that being a twistorial space $(M, \overline{g}, \overline{J})$ has a canonical complex structure (see \cite{Sal3}) and this is 
not hard to see that it is in fact a homogeneous complex contact manifold. It follows then by the results in \cite{Wolf} that the base quaternionic-K\"ahler manifold 
is in fact symmetric and this ends the proof of theorem 8.1.

$\\$
$\\$
$\\$
\begin{flushright}
Paul-Andi Nagy \\
Institut de Math\'ematiques \\
rue E. Argand 11, 2007 Neuch\^atel \\
Switzerland \\
e-mail : Paul.Nagy@unine.ch
\end{flushright}
\end{document}